\documentclass[11pt]{article}
\usepackage{amsmath,amssymb,latexsym, amsfonts, amscd, amsthm}
\usepackage[small,nohug,heads=littlevee]{diagrams}
\diagramstyle[labelstyle=\scriptstyle]

\theoremstyle{plain}

\newtheorem{theorem}{Theorem}[section]

\newtheorem{proposition}[theorem]{Proposition}

\newtheorem{lemma}[theorem]{Lemma}

\theoremstyle{definition}
\newtheorem{definition}[theorem]{Definition}

\newcommand{\mat}[4]{\left( \begin{array}{cc} {#1} & {#2} \\ {#3} & {#4}
\end{array} \right)}
\newcommand{\matt}[9]{\left( \begin{array}{ccc} {#1} & {#2} & {#3} \\ {#4} & {#5} & {#6} \\ {#7} & {#8} & {#9}
\end{array} \right)}

\def\bdf{\begin{defn}}
\def\edf{\end{defn}}

\usepackage{color}

\begin{document}

\title{Characters of some simple supercuspidal representations on split tori}

\author{Moshe Adrian}

\maketitle

\section{Introduction}

Let $F$ be a nonarchimedean local field of characteristic zero, $\mathfrak{o}$ its ring of integers, $\mathfrak{p}$ the maximal ideal, $p$ the residual characteristic, $q$ the order of the residue field, and $\varpi$ a fixed uniformizer of $F$.  Let $G$ be a connected reductive group defined over $F$.  If $\pi$ is an irreducible admissible representation of $G$, we denote by $\theta_{\pi}$ its distribution character, which is a linear functional on $C_c^{\infty}(G)$, the locally constant, compactly supported functions on $G$.  Harish-Chandra showed that $\theta_{\pi}$ can be represented by a locally constant function on the regular semisimple set of $G$, which we will also denote $\theta_{\pi}$.

Suppose that $\pi$ is a supercuspidal representation.  Much is known about $\theta_{\pi}$.  The first supercuspidal characters were computed by Sally and Shalika in \cite{sallyshalika}, where they investigated the supercuspidal representations of $SL(2,F)$ when $p \neq 2$.  Shimizu calculated the supercuspidal characters of $GL(2,F)$ in \cite{shimizu}, for $p \neq 2$.  Kutzko began a study of the supercuspidal characters of $GL(\ell,F)$, $\ell$ a prime (see \cite{kutzko}), when $\ell \neq p$, and DeBacker computed these characters on elliptic tori (see \cite{debacker}).  Later, Spice calculated the supercuspidal characters of $SL(\ell,F)$, $\ell$ an odd prime (see \cite{spice}), for $\ell \neq p$ and together with Adler, they computed a large class of supercuspidal characters for very general connected reductive groups (see \cite{adlerspice}).

Many times one wants to determine character values of a representation on a particular torus, as this can carry much of the information of the representation.  For example, discrete series representations of real groups are determined by their character values on the compact (mod center) torus.  On the $p$-adic side, it is known (see \cite{adrian}) that the supercuspidal representations of $GL(n,F)$, $n$ prime, are determined by their character values on a specific elliptic torus, for $p > 2n$.

In this paper, we compute the character values of the simple supercuspidal representations (recently discovered by Gross and Reeder) of $SL(2,F)$ and $SL(3,F)$ on the maximal split torus of each group when $p$ is arbitrary.  The character values for $SL(2,F)$ are especially elegant, being a fixed constant times a sum of $q^{\ell(w)}$ over appropriate affine Weyl group elements $w$, where $\ell(w)$ is the length of $w$.  We had hoped that the character values of simple supercuspidals for $SL(3,F)$ and more general reductive groups would be as elegant, but this is unfortunately not the case.  As one might expect, various Gauss-type sums appear for $SL(3,F)$.

Let $T$ denote the maximal split torus of $SL(2)$, $Z$ the center, and $W^a$ the affine Weyl group.  Let $val$ denote valuation.  Using the Frobenius character formula for supercuspidal representations (see \cite{sally}), we will first prove the following theorem.

\begin{theorem}\label{maintheorem}
Let $g = \mat{a}{0}{0}{a^{-1}} \in SL(2,F)$ where $a \in 1 + \mathfrak{p}$, and set $r := val(a - a^{-1})$. Let $\pi$ be a simple supercuspidal representation of $SL(2,F)$.  Then
$$\theta_{\pi}(g) = c_q \left( \displaystyle\sum_{w \in W^{a} : \ell(w) < r} q^{\ell(w)} \right)$$
where
\begin{equation*}
c_q := \left\{
\begin{array}{rl}
\frac{q-1}{2} & \text{if } p \neq 2 \\
q-1 & \text{if } p = 2
\end{array} \right.
\end{equation*}
\end{theorem}

It will be clear from our calculations in later sections that $\theta_{\pi}$ vanishes on $T(F) \setminus Z(F) T(1 + \mathfrak{p})$, so up to the central character, we have computed $\theta_{\pi}$ on all of the split torus.  Moreover, since the central character is given by the data forming the simple supercuspidal representation, we have therefore computed all of $\theta_{\pi}$ on the split torus.

We will also compute $\theta_{\pi}(g)$ for simple supercuspidal representations $\pi$ of $SL(3,F)$, where $g$ is in the maximal split torus of $SL(3,F)$.  Specifically, we will compute the character values on $T(1 + \mathfrak{p})$, where here $T$ is the maximal split torus of $SL(3)$.  Then by the same reasoning as for $SL(2,F)$, this will be enough to give the character values on the entire maximal split torus of $SL(3,F)$.  However, since the theorem for $SL(3,F)$ is much more complicated, we will defer its statement to a later section.  Moreover, as there are nontrivial Gauss sums in this formula, we will compute them in complete generality in the last section of this paper.

We note that the term $a - a^{-1}$ in Theorem \ref{maintheorem} is, up to a sign, a canonical square root of the Weyl denominator.  In particular, if $D(g)$ is the standard Weyl denominator, then $-D(g) = (a-a^{-1})^2$.

We wish to make another note.  In \cite{sallyshalika}, Sally/Shalika have character values on the split torus of $SL(2,F)$ for an arbitrary supercuspidal representation of $SL(2,F)$ when $p \neq 2$, which we briefly recall.  For any quadratic extension $V = F(\sqrt{\theta})$ of $F$, let $C_{\theta}$ denote the kernel of the norm $N_{V/F}$, and $\mathfrak{p}_{\theta}$ the prime ideal in $V$.  If $V$ is ramified, set $C_{\theta}^{(h)} = (1 + \mathfrak{p}_{\theta}^{2h+1}) \cap C_{\theta}, h \geq 0$.  If $\psi \in \hat{C}_{\theta}$, denote the conductor of $\psi$ by cond $\psi$ (this is the largest subgroup in the filtration $\{ C_{\theta}^{(h)} \}$ on which $\psi$ is trivial).  The ramified discrete series are indexed by a nontrivial additive character $\eta$ of $F$ and a nontrivial character $\psi \in \hat{C}_{\theta}$, where $V$ is ramified.  The corresponding representation is denoted $\Pi(\eta, \psi, V)$. If $g = \mat{a}{0}{0}{a^{-1}} \in SL(2,F)$, $a \in 1 + \mathfrak{p}$, and $\Pi(\eta, \psi, V)$ is a ramified discrete series, then $$\theta_{\Pi(\eta, \psi, V)}(g) = \frac{1}{|a-a^{-1}|} - \frac{1}{2} q^h \left( \frac{q+1}{q} \right)$$ where $V = F(\sqrt{\theta})$ and cond $\psi = C_{\theta}^{(h)}, h \geq 1$.

Since these character values are in particular valid for the simple supercuspidal representations when $p \neq 2$ (since simple supercuspidal representations are, in particular, ramified discrete series), we may compare their character values to ours.  After some calculation, one can see that their character values agree with ours in the case of simple supercuspidal representations on the split torus.

We now briefly present an outline of the paper.  In section \ref{background}, we define simple supercuspidal representation and we present the relevant background theory that we use to compute the character values.  In section \ref{sl2section}, we compute the character formula for $SL(2,F)$.  In section \ref{sl3section}, we compute the character formula for $SL(3,F)$.  In particular, the formula contains various Gauss-type sums.  In section \ref{morecalculations}, we compute these Gauss sums.

Acknowledgements:
This paper has benefited from conversations with Gordan Savin, Aaron Wood, Chris Kocs, and Loren Spice.

\section{Background}\label{background}

Let us recall some basic definitions.  Let $G$ be a split, simply connected, almost simple, connected reductive group, and $T$ a maximal $F$-split torus in $G$.  Associated to $T$ we have the set of roots $\Phi$ of $T$ in $G$, an apartment, together with a set of affine roots $\Psi$, and an affine Weyl group $W^a$.  We also have a canonical length function $\ell(w)$ on $W^a$.  Fix a Chevalley basis in the Lie algebra of $G$.  To each $\psi \in \Psi$ we have an associated affine root group $U_{\psi}$.   Fix an alcove $C$ in the apartment with corresponding simple and positive affine roots $\Pi \subset \Psi^+$.  Let $T(\mathfrak{o})$ be the maximal compact subgroup of $T(F)$.  Let $T(1+\mathfrak{p}) := < t \in T(\mathfrak{o}) : \lambda(t) \in 1 + \mathfrak{p} \ \forall \lambda \in X^*(T)>$, where $X^*(T)$ is the character lattice of $T$.  Let $I = <T(\mathfrak{o}), U_{\psi} : \psi \in \Psi^+>$ denote the corresponding Iwahori subgroup, and $I_+ = <T(1 + \mathfrak{p}), U_{\psi} : \psi \in \Psi^+>$ its pro-unipotent radical.  Set $I_{++} := <T(1 + \mathfrak{p}), U_{\psi} : \psi \in \Psi^+ \setminus \Pi>$.  We set $H := Z(F) I_+$, where $Z$ is the center of $G$.  Let $N$ denote the normalizer of $T(F)$ in $G(F)$.

\begin{lemma}(see \cite{grossreeder})
The subgroup $I_{++}$ is normal in $I_+$, with quotient $$I_+ / I_{++} \cong \displaystyle\bigoplus_{\psi \in \Pi} U_{\psi} / U_{\psi+1}$$ as $T(\mathfrak{o})$-modules.
\end{lemma}

\begin{definition}(see \cite{grossreeder})
A character $\chi : H \rightarrow \mathbb{C}^*$ is called \emph{affine generic} if

(i) $\chi$ is trivial on $I_{++}$ and

(ii) $\chi$ is nontrivial on $U_{\psi}$ for every $\psi \in \Pi$.
\end{definition}

\begin{theorem} (see \cite{grossreeder})
Let $\chi : H \rightarrow \mathbb{C}^*$ be an affine generic character. Then $cInd_{H}^{G(F)} \chi$ is an irreducible supercuspidal representation, called a \emph{simple supercuspidal representation}, where $cInd$ denotes compact induction.
\end{theorem}

Now suppose that $\pi$ is an irreducible smooth supercuspidal representation of $G(F)$.  Let $K$ be an open, compact subgroup of $G(F)$, and suppose that $\sigma$ is an irreducible representation of $K$ such that $$\pi = cInd_K^{G(F)} \sigma.$$   Let $\chi_{\sigma}$ denote the distribution character of $\sigma$.  The following is the Frobenius formula for the induced character $\theta_{\pi}$.

\begin{theorem}\label{sally} (see \cite{sally})
Let $g$ be a regular element of $G(F)$.  Then $$\theta_{\pi}(g) = \displaystyle\sum_{x \in K \backslash G(F) / K} \displaystyle\sum_{y \in K \backslash K x K} \dot{\chi}_{\sigma}(ygy^{-1})$$ where

\begin{equation*}
\dot{\chi}_{\sigma}(k) = \left\{
\begin{array}{rl}
\chi_{\sigma}(k) & \text{if } k \in K \\
0 & \text{if } k \in G(F) \setminus K
\end{array} \right.
\end{equation*}
\end{theorem}

There is an integral version of this formula as well (see \cite{sally}), and these formulas are a main tool in computing characters of supercuspidal representations.  We will use this theorem to compute the character values of the simple supercuspidal representations of $SL(2,F)$ and $SL(3,F)$ on their maximal split tori.

We first show that the formula in Theorem \ref{sally} simplifies considerably in our situation.  Let us first recall the following basic theory about double coset decompositions.  If $G$ is a connected reductive group and $K$ is a compact open subgroup of $G$, let us choose a set of representatives $\{t_{\alpha} \}$ for the double cosets of $K \backslash G / K$.  Then $K t_{\alpha} K$ is the disjoint union of the cosets $K t_{\alpha} s_1, K t_{\alpha} s_2, ..., K t_{\alpha} s_m$, where $s_1, s_2, ..., s_m$ is a set of representatives of $K / (K \cap t_{\alpha}^{-1} K t_{\alpha})$.  We will use this fact repeatedly in this paper.  We can now state our reduction formula.

\begin{proposition}\label{reduction}
Let $G$ be simply connected, and $\chi$ an affine generic character of $H$.  Set $\pi := cInd_{H}^{G(F)} \chi$.  Then $$\theta_{\pi}(g) = |T(\mathfrak{o}) / Z(F) T(1 + \mathfrak{p})| \displaystyle\sum_{x \in W^a} \displaystyle\sum_{y \in H \backslash H x H} \dot{\chi}_{\sigma}(ygy^{-1}),$$
where the outer sum is meant to be taken over any set of representatives $x$ in $W^a$.
\end{proposition}

\proof
Recall the affine Bruhat decomposition $I \backslash G / I \leftrightarrow W^a \cong N / T(\mathfrak{o})$.  As $I / I_+ \cong T(\mathfrak{o}) / T(1 + \mathfrak{p})$, the affine Bruhat decomposition descends to $H \backslash G / H \leftrightarrow N / Z(F) T(1 + \mathfrak{p})$ (we are using here that $G$ is simply connected, so that $Z(F) = Z(\mathfrak{o})$ and therefore $Z(\mathfrak{o}) I_+ = H \subset I$ and $I / H = T(\mathfrak{o}) / Z(F) T(1 + \mathfrak{p})$), and we have the short exact sequence $$1 \rightarrow T(\mathfrak{o}) / Z(F) T(1 + \mathfrak{p}) \rightarrow N / Z(F) T(1 + \mathfrak{p}) \rightarrow W^a = N / T(\mathfrak{o}) \rightarrow 1$$
Therefore, $$\theta_{\pi}(g) = \displaystyle\sum_{x \in N / Z(F) T(1 + \mathfrak{p})} \displaystyle\sum_{y \in H \backslash H x H} \dot{\chi}_{\sigma}(ygy^{-1})$$
Write $$\sigma(x) := \displaystyle\sum_{y \in H \backslash H x H} \dot{\chi}_{\sigma}(ygy^{-1})$$ for $x \in N$.  We claim that $\sigma(x) = \sigma(xt) \ \forall t \in T(\mathfrak{o})$.  Write $HxH = \displaystyle\cup_{i=1}^m Hx x_i$.  Recall that $x_i$ are representatives of $H / (H \cap x^{-1} Hx)$.  As $x \in N$, $x^{-1} H x = Z(F) I'_+$, where $I'_+$ is the pro-unipotent radical of another Iwahori subgroup $I'$.  It is then easy to see that $H / (H \cap x^{-1} Hx)$ is a direct sum of spaces of the form $U_{\gamma} / U_{\gamma + n}$, where $\gamma \in \Phi^+$ or $\gamma = \gamma' + 1$ where $\gamma' \in \Phi^-$, and where $n$ is a non-negative integer.  Now, $$\sigma(xt) = \displaystyle\sum_{y \in H \backslash H xt H} \dot{\chi}_{\sigma}(ygy^{-1}).$$  Write $HxtH = \displaystyle\cup_{j} Hxt y_j$.  $y_j$ are representatives of $H / (H \cap (xt)^{-1} H xt)$.  Since $t \in T(\mathfrak{o})$, we have that $t^{-1} x^{-1} Hxt = x^{-1} Hx$.  Therefore, in particular, $HxtH = \displaystyle\cup_{i=1}^m Hxt x_i$.  Then $$\sigma(xt) = \displaystyle\sum_{i=1}^m \dot{\chi}_{\sigma}(xtx_i gx_i^{-1} t^{-1} x^{-1}) = $$ $$\displaystyle\sum_{i = 1}^m \dot{\chi}_{\sigma}(xtx_i t^{-1} t g t^{-1} tx_i^{-1} t^{-1} x^{-1}) = \displaystyle\sum_{i=1}^m \dot{\chi}_{\sigma}(xtx_i t^{-1}  g  tx_i^{-1} t^{-1} x^{-1}).$$  Since this sum is over spaces of the form $U_{\gamma} / U_{\gamma+n}$ as noted above, and since conjugation by an element $t \in T(\mathfrak{o})$ preserves $U_{\gamma} / U_{\gamma+n}$, we get $$\displaystyle\sum_{i=1}^m \dot{\chi}_{\sigma}(xtx_i t^{-1}  g  tx_i^{-1} t^{-1} x^{-1}) = \displaystyle\sum_{i=1}^m \dot{\chi}_{\sigma}(xx_i   g  x_i^{-1} x^{-1}) = \sigma(x).$$

Therefore, since $\sigma$ is constant along fibers of the above exact sequence, we have that
$$\theta_{\pi}(g) = |T(\mathfrak{o}) / Z(F) T(1 + \mathfrak{p})| \displaystyle\sum_{x \in W^a} \displaystyle\sum_{y \in H \backslash H x H} \dot{\chi}_{\sigma}(ygy^{-1})$$
\qed

Therefore, we only need to compute $\sigma(x)$ as $x$ varies over a set of representatives of $W^a$.  Moreover, since we are only computing $\theta_{\pi}$ on $T(1+\mathfrak{p})$, we only need the data of the affine generic character $\chi$ on $I_+$ (and not on all of $ZI_+$), since if $g \in T(1+\mathfrak{p})$, then the terms $ygy^{-1}$ in $\theta_{\pi}(g)$ will always live in $I_+$.

\section{The character formula for $SL(2,F)$}\label{sl2section}

In this section we prove Theorem \ref{maintheorem}.  We prove this theorem with a case by case investigation.  We compute the inner sum $\sigma(x)$ in Proposition \ref{reduction} for any set of representatives of $W^a$ by decomposing $HxH$ into a union of left cosets, as in the paragraph that immediately precedes proposition \ref{reduction}.  Afterwards, we sum everything up to get $\theta_{\pi}$.

We fix a Haar measure on $F$ such that $\mathfrak{o}$ has volume $1$, and we use the abbreviation $vol$ to denote volume.  Fix an element $g = \mat{a}{0}{0}{a^{-1}} \in T(1+\mathfrak{p})$, and set $r := val(a - a^{-1})$

\begin{proposition}\label{innersum1}
Let $x = \mat{b}{0}{0}{b^{-1}}$.  Then

\begin{equation*}
\displaystyle\sum_{y \in H \backslash H x H} \dot{\chi}(ygy^{-1}) = \left\{
\begin{array}{rl}
vol(\mathfrak{p}^r)^{-1} vol(\mathfrak{p}^{-2n+r}) & \text{if } val(b) = n \geq 0 \ \mathrm{and} \ 2n < r \\
vol(\mathfrak{p}^r)^{-1} vol(\mathfrak{p}^{2n+r}) & \text{if } val(b) = n < 0 \ \mathrm{and} \ -2n < r \\
0 & \text{otherwise}
\end{array} \right.
\end{equation*}
\end{proposition}

\proof
We first rewrite the double coset $H x H$ as a finite union of single right cosets.  Suppose $val(b) = n > 0.$  Since $$H \cap x^{-1} H x = Z(F) \mat{1 + \mathfrak{p}}{\mathfrak{o}}{\mathfrak{p}^{2n+1}}{1 + \mathfrak{p}},$$ we can obtain an explicit disjoint union $$H = \displaystyle\bigcup_{z \in \mathfrak{p} / \mathfrak{p}^{2n+1}} Z(F) \mat{1 + \mathfrak{p}}{\mathfrak{o}}{\mathfrak{p}^{2n+1}}{1 + \mathfrak{p}} \mat{1}{0}{z}{1},$$  where $$\mat{1 + \mathfrak{p}}{\mathfrak{o}}{\mathfrak{p}^{2n+1}}{1 + \mathfrak{p}} := \left\{ \mat{x_1}{x_2}{x_3}{x_4} \in SL(2,F): x_1,x_4 \in 1 + \mathfrak{p}, x_2 \in \mathfrak{o}, x_3 \in \mathfrak{p}^{2n+1} \right\} $$ (we will use this last type of notation throughout the paper). Therefore, we have a disjoint union $$H x H = \displaystyle\bigcup_{z \in \mathfrak{p} / \mathfrak{p}^{2n+1}} H x \mat{1}{0}{z}{1}$$

Now suppose that $val(b) = n < 0.$  Then similarly, we get $$H \cap x^{-1} H x = Z \mat{1 + \mathfrak{p}}{\mathfrak{p}^{-2n}}{\mathfrak{p}}{1 + \mathfrak{p}},$$ $$H = \displaystyle\bigcup_{z \in \mathfrak{o} / \mathfrak{p}^{-2n}} Z \mat{1 + \mathfrak{p}}{\mathfrak{p}^{-2n}}{\mathfrak{p}}{1 + \mathfrak{p}} \mat{1}{z}{0}{1}.$$  Therefore, we have a disjoint union $$H x H = \displaystyle\bigcup_{z \in \mathfrak{o} / \mathfrak{p}^{-2n}} H x \mat{1}{z}{0}{1}$$

Finally, if $val(b) = 0$, then we get $H \cap x H x^{-1} = H$.  Therefore, $H x H = H x$.

Let us return to the case $val(b) = n > 0$.  Suppose $y \in H x H$.  We need to check when $y g y^{-1} \in H$ since $\dot{\chi}$ vanishes outside $H$. Using our above double coset decomposition, write $y = i x \tilde{z}$, where $\tilde{z}$ is of the form $\mat{1}{0}{z}{1}$ for some $z \in \mathfrak{p} / \mathfrak{p}^{2n+1}$ and for some $i \in H$.  Then $y g y^{-1} \in H \Leftrightarrow x \tilde{z} g \tilde{z}^{-1} x^{-1} \in H$.  Moreover,
$$x \tilde{z} g \tilde{z}^{-1} x^{-1} = \mat{a}{0}{b^{-2} z(a-a^{-1})}{a^{-1}}.$$
Notice that $a \in \pm (1 + \mathfrak{p})$ is forced upon us here in order to have $x \tilde{z} g \tilde{z}^{-1} x^{-1} \in H$.  (We note that the condition $a \in \pm (1 + \mathfrak{p})$ continues to be forced upon us, for the same reason, when you compute the terms $ygy^{-1}$ that appear in $\theta_{\pi}(g)$ for any other representative $x$ of any element of the affine Weyl group,as simple computations will show.  This shows, therefore, that $\theta_{\pi}$ vanishes on $T(F) \setminus Z(F) T(1 + \mathfrak{p})$).

Write $a - a^{-1} = \varpi^r u$ for some unit $u$.  Absorbing all units into the $z$ term, we may write $b^{-2} z(a - a^{-1}) = \varpi^{-2n} \varpi^r z'$ for some $z' \in \mathfrak{p} / \mathfrak{p}^{2n+1}$.

Recall that we are only interested in $\chi|_{I_+}$.  We will abuse notation and write $\chi$ for $\chi|_{I_+}$.  Now, write $\chi$ on $I_+$ as $$\chi : I_+ \rightarrow \mathbb{C}^*$$ $$\mat{d_{11}}{d_{12}}{d_{21}}{d_{22}} \mapsto \chi_1(d_{12}) \chi_2(d_{21}) $$ where $\chi_1$ is a level $1$ character of $\mathfrak{o}$ and where $\chi_2(d_{21}) = \chi_2'(\frac{1}{\varpi} d_{21})$, where $\chi_2'$ is a level $1$ character of $\mathfrak{o}$ (a character of $\mathfrak{o}$ is said to be level $1$ if it is trivial on $\mathfrak{p}$, but nontrivial on $\mathfrak{o}$).  Set $\dot{\chi}_1(z) := \chi_1(z) \ \forall z \in \mathfrak{o}$ and $\dot{\chi}_1(z) = 0 \ \forall z \in F \setminus \mathfrak{o}$.  Moreover, set $\dot{\chi}_2(z) := \chi_2'(z) \ \forall z \in \mathfrak{o}$ and $\dot{\chi}_2(z) = 0 \ \forall z \in F \setminus \mathfrak{o}$.

Therefore,
$$\displaystyle\sum_{y \in H \backslash H x H} \dot{\chi}(ygy^{-1}) = \displaystyle\sum_{z' \in \mathfrak{p} / \mathfrak{p}^{2n+1}} \dot{\chi}_2(\varpi^{-2n+r-1} z')$$ Making a change of variables, we get $$\displaystyle\sum_{z' \in \mathfrak{p} / \mathfrak{p}^{2n+1}} \dot{\chi}_2(\varpi^{-2n+r-1} z') = \displaystyle\sum_{z'' \in \mathfrak{p}^{-2n+r} / \mathfrak{p}^{r}} \dot{\chi}_2(z'') = vol(\mathfrak{p}^r)^{-1} \int_{\mathfrak{p}^{-2n+r} \cap \mathfrak{o}} \dot{\chi}_2(z'') d(z'')$$ since $\dot{\chi}$ vanishes outside $H$.  If $\mathfrak{p}^{-2n+r} \supseteq \mathfrak{o}$, then this integral vanishes since the integral of a nontrivial character over a group vanishes.  However, if $\mathfrak{p}^{-2n+r} \varsubsetneq \mathfrak{o}$, which is precisely the condition that $2n < r$, then $$vol(\mathfrak{p}^r)^{-1} \int_{\mathfrak{p}^{-2n+r} \cap \mathfrak{o}} \dot{\chi}_2(z'') d(z'') = vol(\mathfrak{p}^r)^{-1} \int_{\mathfrak{p}^{-2n+r}} d(z'') = vol(\mathfrak{p}^r)^{-1} vol(\mathfrak{p}^{-2n+r})$$ since $\dot{\chi}_2$ is trivial on $\mathfrak{p}$.

Now consider the case $val(b) = n < 0$.  Suppose $y \in H x H$.  By our above double coset decomposition, write $y = i x \tilde{z}$, where $\tilde{z}$ is of the form $\mat{1}{z}{0}{1}$ for some $z \in \mathfrak{o} / \mathfrak{p}^{-2n}$ and for some $i \in H$.  Moreover, $x \tilde{z} g \tilde{z}^{-1} x^{-1} = \mat{a}{b^2 z(a^{-1} - a)}{0}{a^{-1}}.$  Therefore,
$$\displaystyle\sum_{y \in H \backslash H x H} \dot{\chi}(ygy^{-1}) = \displaystyle\sum_{z \in \mathfrak{o} / \mathfrak{p}^{-2n}} \dot{\chi}_1( b^2 z (a^{-1} - a))$$  We rewrite $b^{2} z(a^{-1} - a) = \varpi^{2n} p^r z'$, where $z' \in \mathfrak{o} / \mathfrak{p}^{-2n}$.
Again, after a change of variables, we get $$\displaystyle\sum_{z \in \mathfrak{o} / \mathfrak{p}^{-2n}} \dot{\chi}_1(b^2 z (a^{-1} - a)) = \displaystyle\sum_{z' \in \mathfrak{o} / \mathfrak{p}^{-2n}} \dot{\chi}_1(\varpi^{2n} p^r z' ) = $$ $$ \displaystyle\sum_{z'' \in \mathfrak{p}^{2n+r} / \mathfrak{p}^{r}} \dot{\chi}_1(z'') = vol(\mathfrak{p}^r)^{-1} \int_{\mathfrak{p}^{2n+r} \cap \mathfrak{o}} \dot{\chi}_1(z'') d(z'')$$  If $\mathfrak{p}^{2n+r} \supseteq \mathfrak{o}$, then again this integral vanishes.  However, if $\mathfrak{p}^{2n+r} \varsubsetneq \mathfrak{o}$, which is precisely the condition that $-2n < r$ then $$vol(\mathfrak{p}^r)^{-1} \int_{\mathfrak{p}^{2n+r} \cap \mathfrak{o}} \dot{\chi}_1(z'') d(z'') = vol(\mathfrak{p}^r)^{-1} \int_{\mathfrak{p}^{2n+r}} d(z'') = vol(\mathfrak{p}^r)^{-1} vol(\mathfrak{p}^{2n+r})$$ since $\dot{\chi}_1$ is trivial on $\mathfrak{p}$.

Now consider the case $val(b) = 0$.  Suppose $y \in H x H$.  Recall that in this case, $H x H = H x$.  Moreover, $x g x^{-1} = g.$  Therefore,
$$\displaystyle\sum_{y \in H \backslash H x H} \dot{\chi}(ygy^{-1}) = \dot{\chi}(g) = 1$$
\qed

\begin{proposition}
Let $x = \mat{0}{c}{-c^{-1}}{0}$.  Then

\begin{equation*}
\displaystyle\sum_{y \in H \backslash H x H} \dot{\chi}(ygy^{-1}) = \left\{
\begin{array}{rl}
vol(\mathfrak{p}^r)^{-1} vol(\mathfrak{p}^{-2n+r-1}) & \text{if } val(c) = n \geq 0 \ \mathrm{and} \ 2n +1 < r \\
vol(\mathfrak{p}^r)^{-1} vol(\mathfrak{p}^{2n+r+1}) & \text{if } val(c) = n < 0 \ \mathrm{and} \ -2n -1 < r \\
0 & \text{otherwise}
\end{array} \right.
\end{equation*}
\end{proposition}

\proof
The proof is completely analogous to that of Proposition \ref{innersum1}.
\qed

Now note that

\begin{equation*}
|T(\mathfrak{o}) / Z(F)T(1+\mathfrak{p})| = |\mathfrak{o}^* / (\pm (1 + \mathfrak{p}))| = \left\{
\begin{array}{rl}
\frac{q-1}{2} & \text{if } p \neq 2 \\
q-1 & \text{if } p = 2
\end{array} \right.
\end{equation*}

(note that if $p = 2$, $1 + \mathfrak{p} = -1 + \mathfrak{p}$).  Thus, our character formula is

\[ \theta_{\pi_{\chi}}(g) = c_q vol(\mathfrak{p}^r)^{-1} \bigg[   \displaystyle\sum_{n \in \mathbb{N}, 0 \leq 2n < r} vol(\mathfrak{p}^{-2n+r}) +   \displaystyle\sum_{n \in \mathbb{N}, 0 < -2n < r} vol(\mathfrak{p}^{2n+r})
\]
\[
+ \displaystyle\sum_{n \in \mathbb{N}, 0 \leq 2n+1 < r} \left( vol(\mathfrak{p}^{-2n+r-1}) \right) + \displaystyle\sum_{n \in \mathbb{N}, 0 < -2n-1 < r} \left( vol(\mathfrak{p}^{2n+r+1}) \right) \bigg]
\]

It is a straightforward calculation to show that if $x = \mat{b}{0}{0}{b^{-1}}$ and $val(b) = k$, then $\ell(x) = |2k|$.  Moreover, if $x = \mat{0}{c}{-c^{-1}}{0}$ and $val(c) = k$, then $\ell(x) = |2k+1|$.

Making the relevant substitutions, and noting that $vol(\mathfrak{p}^d) = q^{-d}$ for $d > 0$ by our choice of measure, one can see that we have proven Theorem \ref{maintheorem}.

\section{The character formula for $SL(3,F)$}\label{sl3section}

In this section we compute the character formula for $SL(3,F)$ on the split maximal torus.  As in the case of $SL(2,F)$, we compute the formula via a case by case investigation.  We compute the inner sum $\sigma(x)$ in Proposition \ref{reduction} for any set of representatives of $W^a$.  Afterwards, we sum everything up to get $\theta_{\pi}$.

Let $g = \matt{\alpha}{0}{0}{0}{\beta}{0}{0}{0}{\gamma} \in T(1+\mathfrak{p})$.  Suppose $\alpha - \beta = \varpi^r u$, $\beta - \gamma = \varpi^s u'$, $\alpha - \gamma = \varpi^t u''$ for some units $u, u', u''$.  Again, we fix a Haar measure on $F$ such that $\mathfrak{o}$ has volume $1$.

Before we state the main theorem, we need to make a few simplifications.  First, notice that $val(\alpha - \gamma) = val( (\alpha - \beta) + (\beta - \gamma) ) \geq inf \{val(\alpha - \beta), val(\beta - \gamma) \}$.  Similarly, $val(\alpha - \beta)  \geq inf \{val(\alpha - \gamma), val(\beta - \gamma) \}$ and $val(\beta - \gamma)  \geq inf \{val(\alpha - \gamma), val(\alpha - \beta) \}$.  One can conclude therefore that either $t \geq r = s$, $s \geq r = t$, or $r \geq s = t$.  Since everything in sight is symmetric, we will assume without loss of generality that $t \geq r = s$, and we will state the main theorem in this case.  One can easily state the analogous results in the two other cases.

So assume $t \geq r = s$.  A complicated impediment is the character values are different in the cases $t = r = s$, $t = r+ 1 = s + 1$, and $t > r + 1 = s + 1$.  We therefore have to state the character formula separately for these three cases.

We need some notation first.  Let $\mathcal{A}_{n_i}$ denote the set of affine Weyl group elements who representatives in $N_G(T)$ are of the form $$\matt{a}{0}{0}{0}{b}{0}{0}{0}{c}, \matt{0}{a}{0}{0}{0}{b}{c}{0}{0}, \matt{0}{0}{a}{b}{0}{0}{0}{c}{0}, $$ $$\matt{a}{0}{0}{0}{0}{b}{0}{c}{0}, \matt{0}{0}{a}{0}{b}{0}{c}{0}{0}, \ \mathrm{and} \ \matt{0}{a}{0}{b}{0}{0}{0}{0}{c},$$ respectively.  For any of the above types of matrices, let $n_{12} = val(a/b), n_{21} = val(b/a), n_{23} = val(b/c), n_{32} = val(c/b), n_{13} = val(a/c), n_{31} = val(c/a)$.  We now define a long list of notation that we need in order to state the main theorem of this section.

Let $\mathcal{B}_1, \mathcal{B}_2, \mathcal{B}_3, \mathcal{B}_4, \mathcal{B}_5, \mathcal{B}_6$ be the inequality conditions $\{n_{21} < r, -n_{31} < t, n_{32} < s \}, \{n_{13} < r-1, n_{21} < s, -n_{23} < t + 1 \}, \{-n_{12} < t + 1, n_{13} < s - 1, n_{32} < r \}, \{n_{21} < t, -n_{31} < r, -n_{23} < s + 1 \}, \{n_{13} < t - 1, -n_{12} < s + 1, -n_{23} < r + 1 \},$ and $\{-n_{12} < r + 1, n_{32} < t, -n_{31} < s \}$, respectively.  For example, if $x = \matt{a}{0}{0}{0}{b}{0}{0}{0}{c}$ and $val(b/a) < r , -val(c/a) < t, val(c/b) < s$, then we say that ``$x$ satisfies $\mathcal{B}_1$''.

Let $\mathcal{C}_1, \mathcal{C}_2, \mathcal{C}_3, \mathcal{C}_4, \mathcal{C}_5, \mathcal{C}_6$ be the inequality conditions $\{n_{31} \leq t, -n_{32} \leq s, -n_{21} \leq r \}, \{n_{23} \leq t - 1, -n_{13} \leq r + 1, -n_{21} \leq s \}, \{-n_{32} \leq r, -n_{13} \leq s + 1, n_{12} \leq t - 1 \}, \{n_{23} \leq s - 1, -n_{21} \leq t, n_{31} \leq r \}, \{n_{23} \leq r - 1, n_{12} \leq s - 1, -n_{13} \leq t + 1 \},$ and $\{n_{31} \leq s, n_{12} \leq r -1, -n_{32} \leq t \}$, respectively.

Let $\mathcal{D}_1, \mathcal{D}_2, \mathcal{D}_3, \mathcal{D}_4, \mathcal{D}_5, \mathcal{D}_6$ be the exact same inequality conditions as $\mathcal{C}_1, \mathcal{C}_2, \mathcal{C}_3, \mathcal{C}_4, \mathcal{C}_5, \mathcal{C}_6$, respectively, except that we replace every $\leq$ sign with a $>$ sign, and moreover, within each $\mathcal{C}_i$, replace every comma by the word ``or''.  For example, if $x = \matt{0}{a}{0}{0}{0}{b}{c}{0}{0}$ satisfies at least one of the inequalities $val(b/c) > t - 1, -val(a/c) > r + 1$, or $-val(b/a) > s$, then we say that ``$x$ satisfies $\mathcal{D}_2$''.

Let $\mathcal{E}_1^0, \mathcal{E}_2^0, \mathcal{E}_3^0, \mathcal{E}_4^0, \mathcal{E}_5^0, \mathcal{E}_6^0$ be the inequality conditions $\{n_{21} \geq 0, n_{31} \geq 0, n_{32} \geq 0 \} \cup \{n_{21} \geq 0, n_{31} < 0, n_{32} < 0 \} \cup \{n_{21} < 0, n_{31} < 0, n_{32} \geq 0 \}, \{n_{13} \geq 0, n_{23} \geq 0, n_{21} \geq 0 \} \cup \{n_{13} < 0, n_{23} < 0, n_{21} \geq 0 \} \cup \{n_{13} \geq 0, n_{23} < 0, n_{21} < 0 \}, \{n_{32} < 0, n_{12} < 0, n_{13} \geq 0 \} \cup \{n_{32} \geq 0, n_{12} < 0, n_{13} < 0 \} \cup \{n_{32} \geq 0, n_{12} \geq 0, n_{13} \geq 0 \}, \{n_{31} < 0, n_{21} \geq 0, n_{23} \geq 0 \} \cup \{n_{31} < 0, n_{21} < 0, n_{23} < 0 \} \cup \{n_{31} \geq 0, n_{21} \geq 0, n_{23} < 0 \}, \{n_{23} \geq 0, n_{13} \geq 0, n_{12} < 0 \} \cup \{n_{23} < 0, n_{13} \geq 0, n_{12} \geq 0 \} \cup \{n_{23} < 0, n_{13} < 0, n_{12} < 0 \},$ and $\{n_{12} < 0, n_{32} \geq 0, n_{31} \geq 0 \} \cup \{n_{12} \geq 0, n_{32} \geq 0, n_{31} < 0 \} \cup \{n_{12} < 0, n_{32} < 0, n_{31} < 0 \}$, respectively.  For example, if $x = \matt{a}{0}{0}{0}{b}{0}{0}{0}{c}$ satisfies $val(b/a) \geq 0, val(c/a) \geq 0$, and $val(c/b) \geq 0$, then we say that ``$x$ satisfies $\mathcal{E}_1^0$.  If $x = \matt{a}{0}{0}{0}{b}{0}{0}{0}{c}$ satisfies $val(b/a) \geq 0, val(c/a) < 0$, and $val(c/b) < 0$, then we say that ``$x$ satisfies $\mathcal{E}_1^0$.  If $x = \matt{a}{0}{0}{0}{b}{0}{0}{0}{c}$ satisfies $val(b/a) < 0, val(c/a) < 0$, and $val(c/b) \geq 0$, then we say that ``$x$ satisfies $\mathcal{E}_1^0$.

Let $\mathcal{E}_1^2, \mathcal{E}_2^2, \mathcal{E}_3^2, \mathcal{E}_4^2, \mathcal{E}_5^2, \mathcal{E}_6^2$ be the inequality conditions $\mathcal{E}_1^0 \setminus \{n_{21} \geq 0, n_{31} \geq 0, n_{32} \geq 0 \}$, $\mathcal{E}_2^0 \setminus \{n_{13} \geq 0, n_{23} \geq 0, n_{21} \geq 0 \}$, $\mathcal{E}_3^0 \setminus \{n_{32} \geq 0, n_{12} \geq 0, n_{13} \geq 0 \}$, $\mathcal{E}_4^0 \setminus \{n_{31} < 0, n_{21} < 0, n_{23} < 0 \}$, $\mathcal{E}_5^0 \setminus \{n_{23} < 0, n_{13} < 0, n_{12} < 0 \}$, and $\mathcal{E}_6^0 \setminus \{n_{12} < 0, n_{32} < 0, n_{31} < 0 \}$, respectively.   For example, if $x = \matt{a}{0}{0}{0}{b}{0}{0}{0}{c}$ satisfies $val(b/a) \geq 0, val(c/a) \geq 0$, and $val(c/b) \geq 0$, then $x$ does not satisfy $\mathcal{E}_1^2$.  If $x = \matt{a}{0}{0}{0}{b}{0}{0}{0}{c}$ satisfies $val(b/a) \geq 0, val(c/a) < 0$, and $val(c/b) < 0$, then ``$x$ satisfies $\mathcal{E}_1^2$.  If $x = \matt{a}{0}{0}{0}{b}{0}{0}{0}{c}$ satisfies $val(b/a) < 0, val(c/a) < 0$, and $val(c/b) \geq 0$, then ``$x$ satisfies $\mathcal{E}_1^2$.

Set $\mathcal{E}_1^1 = \mathcal{E}_1^2, \mathcal{E}_2^1 = \mathcal{E}_2^2, \mathcal{E}_3^1 = \mathcal{E}_3^2, \mathcal{E}_4^1 = \mathcal{E}_4^0, \mathcal{E}_5^1 = \mathcal{E}_5^0, \mathcal{E}_6^1 = \mathcal{E}_6^0$.   For example, if $x = \matt{a}{0}{0}{0}{b}{0}{0}{0}{c}$ satisfies $val(b/a) \geq 0, val(c/a) \geq 0$, and $val(c/b) \geq 0$, then $x$ does not satisfy $\mathcal{E}_1^1$, but it does satisfy $\mathcal{E}_1^0$.

Finally, let $\mathcal{F}_i^j$ be the condition that if $x \in \mathcal{A}_{n_i}$, then $x$ satisfies $\mathcal{B}_i, \mathcal{C}_i,$ and $\mathcal{E}_i^j$.  Moreover, let $\mathcal{G}_i^j$ be the condition that if $x \in \mathcal{A}_{n_i}$, then $x$ satisfies $\mathcal{B}_i, \mathcal{D}_i,$ and $\mathcal{E}_i^j$.  For example, if $x = \matt{a}{0}{0}{0}{0}{b}{0}{c}{0}$, and if $x$ satisfies $val(b/a) < t, -val(c/a) < r, -val(b/c) < s + 1$, and if $x$ satisfies $val(b/c) > s - 1$, and if $x$ satisfies $val(c/a) < 0, val(b/a) \geq 0, val(b/c) \geq 0$, then $x$ satisfies $\mathcal{G}_4^j$, for $j = 0,1,$ and $2$.

We need a few more notations.  As in the case of $SL(2,F)$, the character values will be sums of powers of $q$.  Some of these powers will be lengths of certain affine Weyl group elements, as before, but some will not.  Some powers will be ``truncated lengths'' of affine Weyl group elements.  We will not define ``truncated length''.  However, we will do an example in the next section, and the computation is analogous in every other case.  If an affine Weyl group element $x$ is of type $\mathcal{G}_i^0, \mathcal{G}_i^1$, or $\mathcal{G}_i^2$, then $x$ will contribute the term $q^{\ell'
(x)}$ to $\theta_{\pi}(g)$, where $\ell'(x)$ denotes the ``truncated length'' of $x$.  There will also be two types of Gauss sums that appear in the character formula. If an affine Weyl group element $x$ is of a certain type to be discussed later, then a Gauss sum corresponding to this element will contribute to the character, and we will denote this Gauss sum by either $\Gamma(x)$ or $\Xi(x)$.  We will do examples of how these Gauss sums arise in the next section.

Let $\mathcal{H}_i$ be the condition that if $x \in \mathcal{A}_{n_i}$, then $x$ does not satisfy $\mathcal{E}_i^0$.
Let $\mathcal{E}_i^3$ be the inequality condition $\mathcal{E}_i^0 \setminus \mathcal{E}_i^2$.  Let $\mathcal{J}_i$ be the condition that if $x \in \mathcal{A}_{n_i}$, then $x$ satisfies $\mathcal{E}_i^3$.  Write $\Upsilon := |T(\mathfrak{o}) / Z(F) T(1 + \mathfrak{p})|$.  Finally, in the statement of the main theorem, when we write a summation over $x \in W^a$, we mean that we are summing over any set of representatives of the elements in $W^a$.  Our main theorem for $SL(3,F)$ is the following.

\begin{theorem}
\begin{equation*}
\frac{\theta_{\pi}(g)}{\Upsilon} = \left\{
\begin{array}{ll}
  \displaystyle\sum_{\substack{x \in W^a \\ \mathrm{such \ that}\\x \ \mathrm{satisfies} \ \mathcal{F}_i^0 \\ \mathrm{for \ any} \ 1 \leq i \leq 6}} q^{\ell(x)} + \displaystyle\sum_{\substack{x \in W^a \\ \mathrm{such \ that}\\x \ \mathrm{satisfies} \ \mathcal{G}_i^0 \\ \mathrm{for \ any} \ 1 \leq i \leq 6}} q^{\ell'(x)} + \displaystyle\sum_{\substack{x \in W^a \\ \mathrm{such  \ that} \\ x \ \mathrm{satisfies} \ \mathcal{H}_i \\ \mathrm{for \ any} \ 1 \leq i \leq 6}} \Gamma(x)  & \text{if } t = r \\
 \displaystyle\sum_{\substack{x \in W^a \\ \mathrm{such \ that}\\x \ \mathrm{satisfies} \ \mathcal{F}_i^1 \\ \mathrm{for \ any} \ 1 \leq i \leq 6}} q^{\ell(x)} + \displaystyle\sum_{\substack{x \in W^a \\ \mathrm{such \ that}\\x \ \mathrm{satisfies} \ \mathcal{G}_i^1 \\ \mathrm{for \ any} \ 1 \leq i \leq 6}} q^{\ell'(x)} + \displaystyle\sum_{\substack{x \in W^a \\ \mathrm{such  \ that} \\ x \ \mathrm{satisfies} \ \mathcal{H}_i \\ \mathrm{for \ any} \ 1 \leq i \leq 6}} \Gamma(x) + \displaystyle\sum_{\substack{x \in W^a \\ \mathrm{such \ that}\\x \ \mathrm{satisfies} \ \mathcal{J}_i \\ \mathrm{for \ any} \ 1 \leq i \leq 3}} \Xi(x) & \text{if } t = r + 1 \\
 \displaystyle\sum_{\substack{x \in W^a \\ \mathrm{such \ that}\\x \ \mathrm{satisfies} \ \mathcal{F}_i^2 \\ \mathrm{for \ any} \ 1 \leq i \leq 6}} q^{\ell(x)} + \displaystyle\sum_{\substack{x \in W^a \\ \mathrm{such \ that}\\x \ \mathrm{satisfies} \ \mathcal{G}_i^2 \\ \mathrm{for \ any} \ 1 \leq i \leq 6}} q^{\ell'(x)} + \displaystyle\sum_{\substack{x \in W^a \\ \mathrm{such  \ that} \\ x \ \mathrm{satisfies} \ \mathcal{H}_i \\ \mathrm{for \ any} \ 1 \leq i \leq 6}} \Gamma(x) + \displaystyle\sum_{\substack{x \in W^a \\ \mathrm{such \ that}\\x \ \mathrm{satisfies} \ \mathcal{J}_i \\ \mathrm{for \ any} \ 1 \leq i \leq 6}} \Xi(x) & \text{if } t > r + 1
\end{array} \right.
\end{equation*}
\end{theorem}

We note that the value of term $\Upsilon = |T(\mathfrak{o}) / Z(F) T(1 + \mathfrak{p})|$ depends on whether or not there are cube roots of unity in $F$, which is why we leave this term as is in the above theorem.

\subsection{The case of $x \in \mathcal{A}_{n_1}$}

In this section we compute $\sigma(x)$ when $x \in \mathcal{A}_{n_1}$.  When $x \in \mathcal{A}_{n_2}$ or $x \in \mathcal{A}_{n_3}$, the calculations are similar.  There is a very slight difference, however, when $x \in \mathcal{A}_{n_4}$, $x \in \mathcal{A}_{n_5}$, and $\mathcal{A}_{n_6}$, although these three latter cases are almost completely analogous.  We will address them in the next section.

The main result of this section is the following proposition.  Some of the notation in the proposition will be explained in the proof.

\begin{proposition}\label{firsttheorem}
Let $x = \matt{a}{0}{0}{0}{b}{0}{0}{0}{c}$.  Suppose $n_{21} \geq 0, n_{31} \geq 0, n_{32} \geq 0$. If $r \geq t$, we have

\begin{equation*}
\displaystyle\sum_{y \in H \backslash H x H} \dot{\chi}(ygy^{-1}) = \left\{
\begin{array}{rl}
vol((\mathfrak{o} \cap \mathfrak{p}^{n_{31} - t}) / \mathfrak{p}^{n_{31}}) q^{n_{21} + n_{32}} & \text{if } -n_{21} + r > 0 \ \mathrm{and} \ -n_{32} + s > 0 \ \\
0 & \text{otherwise}
\end{array} \right.
\end{equation*}

If $r < t$, then
\begin{equation*}
\displaystyle\sum_{y \in H \backslash H x H} \dot{\chi}(ygy^{-1}) = \Xi(x)
\end{equation*}

\end{proposition}

\proof
We first rewrite the double coset $H x H$ as a finite union of single right cosets, as in the case of $SL(2,F)$.  Since $$H \cap x^{-1} H x = Z(F) \matt{1 + \mathfrak{p}}{\mathfrak{p}^{n_{21}}}{\mathfrak{p}^{n_{31}}}{\mathfrak{p}}{1 + \mathfrak{p}}{\mathfrak{p}^{n_{32}}}{\mathfrak{p}}{\mathfrak{p}}{1 + \mathfrak{p}},$$ we have a disjoint union $$H x H = \displaystyle\bigcup_{\substack{z_{21} \in \mathfrak{o} / \mathfrak{p}^{n_{21}} \\ z_{32} \in \mathfrak{o} / \mathfrak{p}^{n_{32}} \\ z_{31} \in \mathfrak{o} / \mathfrak{p}^{n_{31}}}} H x \matt{1}{z_{21}}{z_{31}}{0}{1}{z_{32}}{0}{0}{1}$$

Now, if $z = \matt{1}{z_{21}}{z_{31}}{0}{1}{z_{32}}{0}{0}{1}$, then $$xzgz^{-1}x^{-1} = \matt{\alpha}{\frac{a}{b} z_{21} (\beta - \alpha)}{\frac{a}{c}[z_{21} z_{32}(\alpha - \beta) + z_{31}(\gamma - \alpha)]}{0}{\beta}{\frac{b}{c} z_{32} (\gamma - \beta)}{0}{0}{\gamma}$$
Notice that $g \in Z(F) T(1+\mathfrak{p})$ is forced upon us here in order to have $x z g z^{-1} x^{-1} \in H$ (so that $\dot{\chi}$ doesn't vanish on $x z g z^{-1} x^{-1}$).  We note that the condition $g \in Z(F) T(1+\mathfrak{p})$ continues to be forced upon us, for the same reason, when you compute the terms $ygy^{-1}$ that appear in $\theta_{\pi}(g)$ for any other representative $x$ of any element of the affine Weyl group, but we won't include these calculations.
This shows, therefore, that $\theta_{\pi}$ vanishes on $T(F) \setminus Z(F) T(1+\mathfrak{p})$.

Now, write $\chi$ on $I_+$ as $$\chi : I_+ \rightarrow \mathbb{C}^*$$ $$\matt{d_{11}}{d_{12}}{d_{13}}{d_{21}}{d_{22}}{d_{23}}{d_{31}}{d_{32}}{d_{33}} \mapsto \chi_1(d_{12}) \chi_2(d_{23}) \chi_3(d_{31})$$ where $\chi_1, \chi_2$ are level $1$ characters of $\mathfrak{o}$ and where $\chi_3(d_{31}) = \chi_3'(\frac{1}{\varpi} d_{31})$, where $\chi_3'$ is a level $1$ character of $\mathfrak{o}$.

We would like to say that we therefore have
$$\displaystyle\sum_{y \in H \backslash H x H} \dot{\chi}(ygy^{-1}) = \displaystyle\sum_{\substack{z_{21} \in \mathfrak{o} / \mathfrak{p}^{n_{21}} \\ z_{32} \in \mathfrak{o} / \mathfrak{p}^{n_{32}} \\ z_{31} \in \mathfrak{o} / \mathfrak{p}^{n_{31}}}} \dot{\chi}_1\left(\frac{a}{b} z_{21} (\beta - \alpha) \right) \dot{\chi}_2 \left(\frac{b}{c} z_{32} (\gamma - \beta) \right)$$
where $\dot{\chi_i}(z) = \chi_i(z) \ \forall z \in \mathfrak{o}$ and $\dot{\chi_i}(z) = 0 \ \forall z \in F \setminus \mathfrak{o}$, for $i = 1,2$, and where $\dot{\chi}_3(z) = \chi_3'(z) \ \forall z \in \mathfrak{o}$ and $\dot{\chi}_3(z) = 0 \ \forall z \in F \setminus \mathfrak{o}$.

However, we need to take into account the fact that we have the term $\frac{a}{c}[z_{21} z_{32}(\alpha - \beta) + z_{31}(\gamma - \alpha)]$.  $\dot{\chi}$ is zero outside of $Z(F) I_+$, and so we have to take into account the condition that $\frac{a}{c}[z_{21} z_{32}(\alpha - \beta) + z_{31}(\gamma - \alpha)] \in \mathfrak{o}$.

Absorbing all units into the $z_{21}, z_{32}, z_{31}$ terms, and recalling that $\alpha - \beta = \varpi^r u$, $\beta - \gamma = \varpi^s u'$, $\alpha - \gamma = \varpi^t u''$, we therefore wish to understand the condition $\varpi^{-n_{31}+r} z_{21} z_{32} + \varpi^{-n_{31} + t} z_{31} \in \mathfrak{o}$.  We separate this into two cases.

Case 1) Suppose $r \geq t$.  Notice that there might be negative powers of $\varpi$ in $\varpi^{-n_{31}+r} z_{21} z_{32} + \varpi^{-n_{31} + t} z_{31}$ because of the $-n_{31}$ terms.  Fix any $z_{21} \in \mathfrak{o} / \mathfrak{p}^{n_{21}}$ and any $z_{32} \in \mathfrak{o} / \mathfrak{p}^{n_{32}}$.  Then, $\varpi^{-n_{31}+r} z_{21} z_{32}$ can certainly contain negative powers of $\varpi$.  However, we can use $z_{31}$ to cancel out these negative powers of $\varpi$, since $r \geq t$.  In order to force $\varpi^{-n_{31}+r} z_{21} z_{32} + \varpi^{-n_{31} + t} z_{31} \in \mathfrak{o}$, it is evident that all of the negative $\varpi$ power terms in $\varpi^{-n_{31} + t} z_{31}$ are uniquely determined by the negative $\varpi$ power terms in $\varpi^{-n_{31} + r} z_{21} z_{32}$. Moreover, since $r \geq t$, no matter what $z_{21}$ or $z_{32}$ are, we can always find a $z_{31}$ to force $\varpi^{-n_{31} + r} z_{21} z_{32} + \varpi^{-n_{31}+t} z_{31} \in \mathfrak{o}$.  Once we determine the negative $\varpi$ power terms of $\varpi^{-n_{31} + t} z_{31}$ so that the total sum $\varpi^{-n_{31}+r} z_{21} z_{32} + \varpi^{-n_{31} + t} z_{31}$ is in $\mathfrak{o}$, we have  complete leeway in the non-negative $\varpi$ power terms in $\varpi^{-n_{31} + t} z_{31}$.  Therefore, it appears that we have obtained

$$\displaystyle\sum_{y \in H \backslash H x H} \dot{\chi}(ygy^{-1}) = \displaystyle\sum_{\substack{z_{21} \in \mathfrak{o} / \mathfrak{p}^{n_{21}} \\ z_{32} \in \mathfrak{o} / \mathfrak{p}^{n_{32}} \\ z_{31} \in \mathfrak{p}^{n_{31}-t} / \mathfrak{p}^{n_{31}}}} \dot{\chi}_1\left(\frac{a}{b} z_{21} (\beta - \alpha) \right) \dot{\chi}_2 \left(\frac{b}{c} z_{32} (\gamma - \beta) \right)$$

However, the last thing we need to notice is that $z_{31}$ still has to be in $\mathfrak{o}$, by a condition from earlier.  Therefore, the condition $z_{31} \in \mathfrak{p}^{n_{31}-t} / \mathfrak{p}^{n_{31}}$ is not quite correct since it could be the case that $n_{31} < t$.  Therefore, what we really get in then end is

$$\displaystyle\sum_{y \in H \backslash H x H} \dot{\chi}(ygy^{-1}) = \displaystyle\sum_{\substack{z_{21} \in \mathfrak{o} / \mathfrak{p}^{n_{21}} \\ z_{32} \in \mathfrak{o} / \mathfrak{p}^{n_{32}} \\ z_{31} \in (\mathfrak{o} \cap \mathfrak{p}^{n_{31}-t}) / \mathfrak{p}^{n_{31}}}} \dot{\chi}_1\left(\frac{a}{b} z_{21} (\beta - \alpha) \right) \dot{\chi}_2 \left(\frac{b}{c} z_{32} (\gamma - \beta) \right)$$

Absorbing all units into the $z_{21}, z_{32}$ terms, we get
$$vol((\mathfrak{o} \cap \mathfrak{p}^{n_{31} - t}) / \mathfrak{p}^{n_{31}}) \displaystyle\sum_{\substack{z_{21} \in \mathfrak{o} / \mathfrak{p}^{n_{21}} \\ z_{32} \in \mathfrak{o} / \mathfrak{p}^{n_{32}}}} \dot{\chi}_1\left(\varpi^{-n_{21}} z_{21} \varpi^r \right) \dot{\chi}_2 \left(\varpi^{-n_{32}} z_{32} \varpi^{s} \right)$$  Making a change of variables, we get
$$vol((\mathfrak{o} \cap \mathfrak{p}^{n_{31} - t}) / \mathfrak{p}^{n_{31}}) \displaystyle\sum_{\substack{z_{21}' \in \mathfrak{p}^{-n_{21} + r} / \mathfrak{p}^{r} \\ z_{32}' \in \mathfrak{p}^{-n_{32} + s} / \mathfrak{p}^{s}}} \dot{\chi}_1\left(z_{21}' \right) \dot{\chi}_2 \left(z_{32}' \right) = $$
$$vol((\mathfrak{o} \cap \mathfrak{p}^{n_{31} - t}) / \mathfrak{p}^{n_{31}}) \left( vol(\mathfrak{p}^r)^{-1} \displaystyle\int_{\mathfrak{p}^{-n_{21}+r} \cap \mathfrak{o}} \dot{\chi}_1(z_{21}') dz_{21}' \right) \left( vol(\mathfrak{p}^s)^{-1} \displaystyle\int_{\mathfrak{p}^{-n_{32}+s} \cap \mathfrak{o}} \dot{\chi}_2(z_{32}') dz_{32}' \right)$$
Therefore, since the integral of a nontrivial character over a group vanishes, we get

\begin{equation*}
\displaystyle\sum_{y \in H \backslash H x H} \dot{\chi}(ygy^{-1}) = \left\{
\begin{array}{rl}
vol((\mathfrak{o} \cap \mathfrak{p}^{n_{31} - t}) / \mathfrak{p}^{n_{31}}) q^{n_{21} + n_{32}} & \text{if } -n_{21} + r > 0 \ \mathrm{and} \ -n_{32} + s > 0 \ \\
0 & \text{otherwise}
\end{array} \right.
\end{equation*}

which concludes Case 1).

Case 2): Suppose $r < t$.  In this case, if one picks any random $z_{21}$ and $z_{32}$, one might not be able to find a $z_{31}$ that makes $\varpi^{-n_{31}+r} z_{21} z_{32} + \varpi^{-n_{31} + t} z_{31}$ land in $\mathfrak{o}$.  However, if one chooses $z_{21}, z_{32}$ such that $z_{21} z_{32} \in \mathfrak{p}^{t-r}$, then one can find a $z_{31}$ such that $\varpi^{-n_{31}+r} z_{21} z_{32} + \varpi^{-n_{31} + t} z_{31}$ will be in $\mathfrak{o}$.  One can see that $z_{21} z_{32} \in \mathfrak{p}^{t-r}$ is the only additional condition that we need to add to the conditions in Case 1, so we get

$$\displaystyle\sum_{y \in H \backslash H x H} \dot{\chi}(ygy^{-1}) = $$ $$vol((\mathfrak{o} \cap \mathfrak{p}^{n_{31} - t}) / \mathfrak{p}^{n_{31}}) \displaystyle\sum_{\substack{z_{21} \in \mathfrak{o} / \mathfrak{p}^{n_{21}}, z_{32} \in \mathfrak{o} / \mathfrak{p}^{n_{32}} \\ \mathrm{such \ that \ } z_{21} z_{32} \in \mathfrak{p}^{t-r}}} \dot{\chi}_1\left(\frac{a}{b} z_{21} (\beta - \alpha) \right) \dot{\chi}_2 \left(\frac{b}{c} z_{32} (\gamma - \beta) \right)$$

Absorbing all units into the $z_{21}, z_{32}$ terms, we get
\begin{equation}
vol((\mathfrak{o} \cap \mathfrak{p}^{n_{31} - t}) / \mathfrak{p}^{n_{31}}) \displaystyle\sum_{\substack{z_{21} \in \mathfrak{o} / \mathfrak{p}^{n_{21}}, z_{32} \in \mathfrak{o} / \mathfrak{p}^{n_{32}} \\ \mathrm{such \ that \ } z_{21} z_{32} \in \mathfrak{p}^{t-r}}} \dot{\chi}_1\left(\varpi^{-n_{21}} z_{21} \varpi^r \right) \dot{\chi}_2 \left(\varpi^{-n_{32}} z_{32} \varpi^{s} \right) \ \ \label{xi}
\end{equation}
Since this sum is quite complicated, we compute this type of sum in full generality in a later section.  We will instead denote the value of $\displaystyle\sum_{y \in H \backslash H x H} \dot{\chi}(ygy^{-1})$ by $\Xi(x)$ for the $x$ in this proposition, in the case that $r < t$.  Moreover, since we will encounter the analogous type of sum in equation (\ref{xi}) for many other elements $x$, we will merely denote the value of $\displaystyle\sum_{y \in H \backslash H x H} \dot{\chi}(ygy^{-1})$ by $\Xi(x)$ for those $x$ as well.
\qed

We note that in Case 1) above, if $n_{31} \leq t$, we get $vol((\mathfrak{o} \cap \mathfrak{p}^{n_{31} - t}) / \mathfrak{p}^{n_{31}}) q^{n_{21} + n_{32}} = q^{n_{31} + n_{21} + n_{32}}$, which equals $q^{\ell(x)}$, as a simple calculation will show.  This is how the terms of the form $q^{\ell(x)}$ appear in the character formula.  If, on the other hand, $n_{31} > t$, we get $vol((\mathfrak{o} \cap \mathfrak{p}^{n_{31} - t}) / \mathfrak{p}^{n_{31}}) q^{n_{21} + n_{32}} = q^{t + n_{21} + n_{32}}$, which is what we have called $q^{\ell'(x)}$, as $t + n_{21} + n_{32}$ is a ``truncated length'' of $x$.

We now illustrate one more case in the $\mathcal{A}_{n_1}$ setting, which will show how the second type of Gauss sum in the character formula arises. Some of the notation in the next proposition will be explained in the proof.

\begin{proposition}
Let $x = \matt{a}{0}{0}{0}{b}{0}{0}{0}{c}$.  Suppose $n_{21} < 0, n_{31} < 0, n_{32} < 0$. Then

\begin{equation*}
\displaystyle\sum_{y \in H \backslash H x H} \dot{\chi}(ygy^{-1}) = \Gamma(x)
\end{equation*}
\end{proposition}

\proof
We have a disjoint union $$H x H = \displaystyle\bigcup_{\substack{z_{21} \in \mathfrak{p} / \mathfrak{p}^{-n_{21}+1} \\ z_{31} \in \mathfrak{p} / \mathfrak{p}^{-n_{31}+1} \\ z_{32} \in \mathfrak{p} / \mathfrak{p}^{n_{32}+1}}} H x \matt{1}{0}{0}{z_{21}}{1}{0}{z_{31}}{z_{32}}{1}$$

Now, if $z = \matt{1}{0}{0}{z_{21}}{1}{0}{z_{31}}{z_{32}}{1}$, then $$xzgz^{-1}x^{-1} = \matt{\alpha}{0}{0}{\frac{b}{a} z_{21} (\alpha - \beta)}{\beta}{0}{\frac{c}{a} [ z_{32} z_{21} (\gamma - \beta) + z_{31} ( \alpha - \gamma) ]}{\frac{c}{b} z_{32} ( \beta - \gamma)}{\gamma}$$

\begin{equation}
\displaystyle\sum_{y \in H \backslash H x H} \dot{\chi}(ygy^{-1}) = \displaystyle\sum \dot{\chi}(\frac{c}{a} [ z_{32} z_{21} (\gamma - \beta) + z_{31} ( \alpha - \gamma) ] \ \ \label{gamma}
\end{equation}

where the sum is over $z_{21} \in \mathfrak{p} / \mathfrak{p}^{-n_{21}+1}, z_{31} \in \mathfrak{p} / \mathfrak{p}^{-n_{31}+1}, z_{32} \in \mathfrak{p} / \mathfrak{p}^{n_{32}+1} \ \mathrm{such \ that} \ \frac{b}{a} z_{21} (\alpha - \beta) \in \mathfrak{p} \ \mathrm{and} \ \frac{c}{b} z_{32} ( \beta - \gamma) \in \mathfrak{p}$.  Since this type of sum is quite complicated, we compute this type of sum in full generality in section \ref{morecalculations}.  We will instead denote the value of $\displaystyle\sum_{y \in H \backslash H x H} \dot{\chi}(ygy^{-1})$ by $\Gamma(x)$ for the $x$ in this proposition.  Moreover, since we will encounter the type of sum on the right hand side of the above equation (\ref{gamma}) for many other elements $x$, we will merely denote the value of $\displaystyle\sum_{y \in H \backslash H x H} \dot{\chi}(ygy^{-1})$ by $\Gamma(x)$ for those $x$ as well.

\subsection{The case of $x \in \mathcal{A}_{n_4}$}
In this section we will consider the inner sums $\sigma(x)$ when $x \in \mathcal{A}_{n_4}$.  This case is slightly different from that of $\mathcal{A}_{n_1}, \mathcal{A}_{n_2}$, and $\mathcal{A}_{n_3}$.  In particular, the inequalities between $r,s,$ and $t$ that distinguished between Case 1)'s and Case 2)'s in the previous section are now shifted, as we shall show.  This is why we need separate cases in the statement of the main theorem of the distribution character for $SL(3,F)$.  We will show why the shifts occur in this section.  The cases of $x \in \mathcal{A}_{n_5}$ and $x \in \mathcal{A}_{n_6}$ are similar to the case of $x \in \mathcal{A}_{n_4}$.

\begin{proposition}
Let $x = \matt{a}{0}{0}{0}{0}{b}{0}{c}{0}$.  Suppose $n_{31} < 0, n_{21} < 0, n_{23} < 0$. Then if $s \geq t -1$, we have
$$\displaystyle\sum_{y \in H \backslash H x H} \dot{\chi}(ygy^{-1}) =$$

\begin{equation*}
\left\{
\begin{array}{rl}
vol((\mathfrak{p} \cap \mathfrak{p}^{-n_{21} - t + 1}) / \mathfrak{p}^{-n_{21}+1}) q^{-n_{31} - n_{23}-1} & \text{if } n_{23}+s+1 > 0 \ \mathrm{and} \ n_{31}+r > 0 \ \\
0 & \text{otherwise}
\end{array} \right.
\end{equation*}

If $s < t-1$, then
\begin{equation*}
\displaystyle\sum_{y \in H \backslash H x H} \dot{\chi}(ygy^{-1}) = \Xi(x)
\end{equation*}
\end{proposition}

\proof
We have a disjoint union $$H x H = \displaystyle\bigcup_{\substack{z_{21} \in \mathfrak{p} / \mathfrak{p}^{-n_{21}+1} \\ z_{23} \in \mathfrak{p} / \mathfrak{p}^{-n_{23}} \\ z_{31} \in \mathfrak{p} / \mathfrak{p}^{-n_{31}+1}}} H x \matt{1}{0}{0}{z_{31}}{1}{0}{z_{21}}{z_{23}}{1}$$

Now, if $z = \matt{1}{0}{0}{z_{31}}{1}{0}{z_{21}}{z_{23}}{1}$, then $$xzgz^{-1}x^{-1} = \matt{\alpha}{0}{0}{\frac{b}{a}[z_{21}(\alpha - \gamma) + z_{31} z_{23} (\gamma - \beta)]}{\gamma}{\frac{b}{c} z_{23} (\beta - \gamma)}{\frac{c}{a} z_{31} (\alpha - \beta)}{0}{\beta}$$

We have to take into account the condition that $\frac{b}{a}[z_{21}(\alpha - \gamma) + z_{31} z_{23} (\gamma - \beta)] \in \mathfrak{p}$.  We therefore wish to understand the condition $\varpi^{n_{21}+s} z_{31} z_{23} + \varpi^{n_{21} + t} z_{21} \in \mathfrak{p}$.  We separate this into two cases, exactly in the way we did in Proposition \ref{firsttheorem}.  Notice here that $z_{21} \in \mathfrak{p}, z_{23} \in \mathfrak{p}, z_{31} \in \mathfrak{p}$.  Therefore, our first case is going to be

Case 1) Suppose $s \geq t -1$.  Then an analogous computation as we have done before shows that $$\displaystyle\sum_{y \in H \backslash H x H} \dot{\chi}(ygy^{-1}) =$$

\begin{equation*}
\left\{
\begin{array}{rl}
vol((\mathfrak{p} \cap \mathfrak{p}^{-n_{21} - t + 1}) / \mathfrak{p}^{-n_{21}+1}) q^{-n_{31} - n_{23}-1} & \text{if } n_{23}+s+1 > 0 \ \mathrm{and} \ n_{31}+r > 0 \ \\
0 & \text{otherwise}
\end{array} \right.
\end{equation*}

Case 2) Suppose $s < t - 1$.  Then, using notation from the previous section, we have $$\displaystyle\sum_{y \in H \backslash H x H} \dot{\chi}(ygy^{-1}) = \Xi(x)$$
\qed

\section{Calculation of the Gauss sums}\label{morecalculations}

\subsection{The calculation of $\Gamma(x)$}\label{gammasections}

In this section we calculate the terms of the form $\Gamma(x)$ in a much more general context.  In particular, we shall calculate $$\displaystyle\sum_{\substack{x \in \mathfrak{p}^m / \mathfrak{p}^n \\ y \in \mathfrak{p}^k / \mathfrak{p}^{\ell} \\ z \in \mathfrak{p}^i / \mathfrak{p}^j}   } \chi(\varpi^{-a} xy + \varpi^{-b} z),$$ where $\chi$ is a character of $F$ that is zero on $F \setminus \mathfrak{o}$, and is $1$ on $\mathfrak{p}$, except when

i) $m + k < a, i < b$, and $n + \ell \leq a$

ii) $m + k < a, i < b$, and $j \leq b$

The reason that we may omit these cases is that the conditions $j \leq b$ and $n + \ell \leq a$ never occur in the $SL(3,F)$ calculations (as one can check, but we do not include these details here), so we can ignore them.  We won't ignore every possible case where $j \leq b$ or $n + \ell \leq a$, since some of these cases are easy to write down.  We will just ignore the above two special cases.  We split up the calculation of the above sum into various cases.

\

Case 1) Suppose $m+k \geq a, i > b$.  We then have that $\varpi^{-a} xy, \varpi^{-b} z \in \mathfrak{o}$.  Therefore, $\chi(\varpi^{-a} xy + \varpi^{-b} z) = \chi(\varpi^{-a} xy) \chi( \varpi^{-b} z)$, so $$\displaystyle\sum_{\substack{x \in \mathfrak{p}^m / \mathfrak{p}^n \\ y \in \mathfrak{p}^k / \mathfrak{p}^{\ell} \\ z \in \mathfrak{p}^i / \mathfrak{p}^j}   } \chi(\varpi^{-a} xy + \varpi^{-b} z) = \displaystyle\sum_{\substack{x \in \mathfrak{p}^m / \mathfrak{p}^n \\ y \in \mathfrak{p}^k / \mathfrak{p}^{\ell}}} \chi(\varpi^{-a} xy) \displaystyle\sum_{z \in \mathfrak{p}^i / \mathfrak{p}^j } \chi( \varpi^{-b} z)$$
After making a change of variables, we get $$\displaystyle\sum_{z \in \mathfrak{p}^i / \mathfrak{p}^j } \chi( \varpi^{-b} z) = \displaystyle\sum_{z' \in \mathfrak{p}^{i-b} / \mathfrak{p}^{j-b} } \chi( z') $$ $$=vol(\mathfrak{p}^{j-b})^{-1} \int_{\mathfrak{p}^{i-b} \cap \mathfrak{o}} \chi(z') dz' = vol(\mathfrak{p}^{j-b})^{-1} vol(\mathfrak{p}^{i-b}) = q^{j-i}$$ and so $$\displaystyle\sum_{\substack{x \in \mathfrak{p}^m / \mathfrak{p}^n \\ y \in \mathfrak{p}^k / \mathfrak{p}^{\ell} \\ z \in \mathfrak{p}^i / \mathfrak{p}^j} } \chi(\varpi^{-a} xy + \varpi^{-b} z) = q^{j-i} \displaystyle\sum_{\substack{x \in \mathfrak{p}^m / \mathfrak{p}^n \\ y \in \mathfrak{p}^k / \mathfrak{p}^{\ell}}} \chi(\varpi^{-a} xy)$$  We now consider two subcases : $m + k > a$ and $m + k = a$.  Suppose that $m + k > a$.  Then $\varpi^{-a} xy \in \mathfrak{p}$, and therefore $\chi(\varpi^{-a} xy) = 1 \ \forall x \in \mathfrak{p}^m / \mathfrak{p}^n, y \in \mathfrak{p}^k / \mathfrak{p}^{\ell}$ since $\chi$ is trivial on $\mathfrak{p}$.  Therefore, $$\displaystyle\sum_{\substack{x \in \mathfrak{p}^m / \mathfrak{p}^n \\ y \in \mathfrak{p}^k / \mathfrak{p}^{\ell}}} \chi(\varpi^{-a} xy) = \displaystyle\sum_{\substack{x \in \mathfrak{p}^m / \mathfrak{p}^n \\ y \in \mathfrak{p}^k / \mathfrak{p}^{\ell}}} 1 = vol(\mathfrak{p}^m / \mathfrak{p}^n) vol(\mathfrak{p}^k / \mathfrak{p}^{\ell}) = q^{n-m} q^{\ell-k}$$
Now suppose $m + k = a$.  We argue by fixing values of $x$.  Suppose $val(x) = m$.  Then
$$\displaystyle\sum_{y \in \mathfrak{p}^k / \mathfrak{p}^{\ell}} \chi(\varpi^{-a} xy) = 0$$ since $\varpi^{-a} xy$ ranges over all elements of $\mathfrak{o} / \mathfrak{p}^{\ell-a+m}$.  Therefore, the contributions to $$\displaystyle\sum_{\substack{x \in \mathfrak{p}^m / \mathfrak{p}^n \\ y \in \mathfrak{p}^k / \mathfrak{p}^{\ell}}} \chi(\varpi^{-a} xy)$$ from the elements $x$ that have valuation $m$ are all zeroes.  Therefore,
$$\displaystyle\sum_{\substack{x \in \mathfrak{p}^m / \mathfrak{p}^n \\ y \in \mathfrak{p}^k / \mathfrak{p}^{\ell}}} \chi(\varpi^{-a} xy) = \displaystyle\sum_{\substack{x \in \mathfrak{p}^{m+1} / \mathfrak{p}^n \\ y \in \mathfrak{p}^k / \mathfrak{p}^{\ell}}} \chi(\varpi^{-a} xy)$$  But now in this sum, notice that since we have $x \in \mathfrak{p}^{m+1} / \mathfrak{p}^n, y \in \mathfrak{p}^k / \mathfrak{p}^{\ell}$, we conclude that $\varpi^{-a} xy$ is always in $\mathfrak{p}$.  Therefore, as before,
$$\displaystyle\sum_{\substack{x \in \mathfrak{p}^{m+1} / \mathfrak{p}^n \\ y \in \mathfrak{p}^k / \mathfrak{p}^{\ell}}} \chi(\varpi^{-a} xy) = \displaystyle\sum_{\substack{x \in \mathfrak{p}^{m+1} / \mathfrak{p}^n \\ y \in \mathfrak{p}^k / \mathfrak{p}^{\ell}}} 1 = q^{n - (m+1)} q^{\ell-k}$$

Case 2) Suppose $m+k \geq a, i = b$.  This case is mostly analagous to Case 1), except we now have $$\displaystyle\sum_{z \in \mathfrak{p}^i / \mathfrak{p}^j } \chi( \varpi^{-b} z) = 0$$ since $i = b$, and therefore $$\displaystyle\sum_{\substack{x \in \mathfrak{p}^m / \mathfrak{p}^n \\ y \in \mathfrak{p}^k / \mathfrak{p}^{\ell} \\ z \in \mathfrak{p}^i / \mathfrak{p}^j} } \chi(\varpi^{-a} xy + \varpi^{-b} z) = 0$$

Case 3) Suppose $m+k \geq a, i < b$.  We need to understand the condition $\varpi^{-a} xy + \varpi^{-b} z \in \mathfrak{o}$. We first assume that $j > b$.  Thus,
$$\displaystyle\sum_{\substack{x \in \mathfrak{p}^m / \mathfrak{p}^n \\ y \in \mathfrak{p}^k / \mathfrak{p}^{\ell} \\ z \in \mathfrak{p}^i / \mathfrak{p}^j }} \chi(\varpi^{-a} xy + \varpi^{-b} z) = \displaystyle\sum_{\substack{x \in \mathfrak{p}^m / \mathfrak{p}^n \\ y \in \mathfrak{p}^k / \mathfrak{p}^{\ell} \\ z \in \mathfrak{p}^b / \mathfrak{p}^j }} \chi(\varpi^{-a} xy + \varpi^{-b} z)$$  (note that $b < j$, so we can talk about $\mathfrak{p}^b / \mathfrak{p}^j$).  But now note that since $z \in \mathfrak{p}^b$, we have that $\varpi^{-a} xy, \varpi^{-b} z \in \mathfrak{o}$.  Therefore,
$$ \displaystyle\sum_{\substack{x \in \mathfrak{p}^m / \mathfrak{p}^n \\ y \in \mathfrak{p}^k / \mathfrak{p}^{\ell} \\ z \in \mathfrak{p}^b / \mathfrak{p}^j }} \chi(\varpi^{-a} xy + \varpi^{-b} z) =
\displaystyle\sum_{\substack{x \in \mathfrak{p}^m / \mathfrak{p}^n \\ y \in \mathfrak{p}^k / \mathfrak{p}^{\ell}}} \chi(\varpi^{-a} xy) \displaystyle\sum_{ z \in \mathfrak{p}^b / \mathfrak{p}^j } \chi(\varpi^{-b} z) = 0$$
since the integral over a group of a nontrivial character vanishes.

If $j \leq b$, then in order for $\varpi^{-a} xy + \varpi^{-b} z$ to be in $\mathfrak{o}$, we require that $z = 0$, since $z$ is assumed to be in $\mathfrak{p}^i / \mathfrak{p}^j$.  Therefore,
$$\displaystyle\sum_{\substack{x \in \mathfrak{p}^m / \mathfrak{p}^n \\ y \in \mathfrak{p}^k / \mathfrak{p}^{\ell} \\ z \in \mathfrak{p}^i / \mathfrak{p}^j }} \chi(\varpi^{-a} xy + \varpi^{-b} z) = \displaystyle\sum_{\substack{x \in \mathfrak{p}^m / \mathfrak{p}^n \\ y \in \mathfrak{p}^k / \mathfrak{p}^{\ell}}} \chi(\varpi^{-a} xy),$$ which may be rewritten as $$\displaystyle\sum_{\substack{x' \in \mathfrak{p}^{m-a} / \mathfrak{p}^{n-a} \\ y \in \mathfrak{p}^k / \mathfrak{p}^{\ell}}} \chi(x'y)$$ after a change of variables.  This type of sum will be handled in Case 4), which we now present.

Case 4) Suppose $m+k < a, i > b$.  We consider two subcases.  We therefore get $$\displaystyle\sum_{\substack{x \in \mathfrak{p}^m / \mathfrak{p}^n \\ y \in \mathfrak{p}^k / \mathfrak{p}^{\ell} \\ z \in \mathfrak{p}^i / \mathfrak{p}^j }} \chi(\varpi^{-a} xy + \varpi^{-b} z) = \displaystyle\sum_{\substack{x \in \mathfrak{p}^m / \mathfrak{p}^n, y \in \mathfrak{p}^k / \mathfrak{p}^{\ell}, z \in \mathfrak{p}^i / \mathfrak{p}^j \\ \mathrm{such \ that \ } xy \in \mathfrak{p}^a }} \chi(\varpi^{-a} xy + \varpi^{-b} z)$$     since if $xy \notin \mathfrak{p}^a$, then $\varpi^{-a} xy + \varpi^{-b} z \notin \mathfrak{o}$.
Therefore, we get $$\displaystyle\sum_{\substack{x \in \mathfrak{p}^m / \mathfrak{p}^n \\ y \in \mathfrak{p}^k / \mathfrak{p}^{\ell} \\ z \in \mathfrak{p}^i / \mathfrak{p}^j} } \chi(\varpi^{-a} xy + \varpi^{-b} z) = \displaystyle\sum_{\substack{x \in \mathfrak{p}^m / \mathfrak{p}^n, y \in \mathfrak{p}^k / \mathfrak{p}^{\ell}, z \in \mathfrak{p}^i / \mathfrak{p}^j \\ \mathrm{such \ that \ } xy \in \mathfrak{p}^a }} \chi(\varpi^{-a} xy + \varpi^{-b} z) =$$ $$ \displaystyle\sum_{\substack{x \in \mathfrak{p}^m / \mathfrak{p}^n, y \in \mathfrak{p}^k / \mathfrak{p}^{\ell} \\ \mathrm{such \ that \ } xy \in \mathfrak{p}^a}} \chi(\varpi^{-a} xy) \displaystyle\sum_{ z \in \mathfrak{p}^i / \mathfrak{p}^j } \chi(\varpi^{-b} z) = q^{j-i} \displaystyle\sum_{\substack{x \in \mathfrak{p}^m / \mathfrak{p}^n, y \in \mathfrak{p}^k / \mathfrak{p}^{\ell} \\ \mathrm{such \ that \ } xy \in \mathfrak{p}^a}} \chi(\varpi^{-a} xy)$$ We make a change of variables  $$\displaystyle\sum_{\substack{x \in \mathfrak{p}^m / \mathfrak{p}^n, y \in \mathfrak{p}^k / \mathfrak{p}^{\ell} \\ \mathrm{such \ that \ } xy \in \mathfrak{p}^a}} \chi(\varpi^{-a} xy) = \displaystyle\sum_{\substack{x' \in \mathfrak{p}^{m-a} / \mathfrak{p}^{n-a}, y \in \mathfrak{p}^k / \mathfrak{p}^{\ell} \\ \mathrm{such \ that \ } x'y \in \mathfrak{o}}} \chi(x'y)$$

Since $\chi$ vanishes outside $\mathfrak{o}$ $$\displaystyle\sum_{\substack{x' \in \mathfrak{p}^{m-a} / \mathfrak{p}^{n-a}, y \in \mathfrak{p}^k / \mathfrak{p}^{\ell} \\ \mathrm{such \ that \ } x'y \in \mathfrak{o}}} \chi(x'y) = \displaystyle\sum_{\substack{x' \in \mathfrak{p}^{m-a} / \mathfrak{p}^{n-a} \\ y \in \mathfrak{p}^k / \mathfrak{p}^{\ell}}} \chi(x'y)$$

After reindexing and relabeling, we are now interested in computing the following type of sum $$\displaystyle\sum_{\substack{x \in \mathfrak{p}^{m} / \mathfrak{p}^{n} \\ y \in \mathfrak{p}^k / \mathfrak{p}^{\ell}}} \chi(xy)$$
We split this into two cases:

Case i) Suppose $n + \ell \geq  -2$.  Recall that $\chi$ vanishes outside $\mathfrak{o}$, so we are interested in when $xy \in \mathfrak{o}$.  We will separate out the $x$ terms from the sum that have no chance of multiplying with an element of $y$ to land in $\mathfrak{o}$ unless $y$ is zero.  We consider two subcases.  Suppose $1 - \ell > m$.   We write $$\displaystyle\sum_{\substack{x \in \mathfrak{p}^{m} / \mathfrak{p}^{n} \\ y \in \mathfrak{p}^k / \mathfrak{p}^{\ell}}} \chi(xy) = \displaystyle\sum_{\substack{x \in \mathfrak{p}^{m} / \mathfrak{p}^{n} - \mathfrak{p}^{1-\ell} / \mathfrak{p}^n \\ y \in \mathfrak{p}^k / \mathfrak{p}^{\ell}}} \chi(xy) + \displaystyle\sum_{\substack{x \in \mathfrak{p}^{1-\ell} / \mathfrak{p}^{n} \\ y \in \mathfrak{p}^k / \mathfrak{p}^{\ell}}} \chi(xy)$$  (note that we can write $x \in \mathfrak{p}^{m} / \mathfrak{p}^{n} - \mathfrak{p}^{1-\ell} / \mathfrak{p}^n$ since we assumed $1 - \ell \geq m$).  Let the first sum be denoted $A$ and the second sum be denoted $B$. For the first sum $A$, $xy$ can never be in $\mathfrak{o}$ unless $y = 0$.  Therefore, we get $$A = \displaystyle\sum_{\substack{x \in \mathfrak{p}^{m} / \mathfrak{p}^{n} - \mathfrak{p}^{1-\ell} / \mathfrak{p}^n \\ y = 0}} \chi(xy) =  \displaystyle\sum_{x \in \mathfrak{p}^{m} / \mathfrak{p}^{n} - \mathfrak{p}^{1-\ell} / \mathfrak{p}^n} 1 = q^{n-m} - q^{n - (1 - \ell)}.$$   For $B$, we split up the sum as

$$\displaystyle\sum_{\substack{x \in \mathfrak{p}^{1-\ell} / \mathfrak{p}^{n} \\ y \in \mathfrak{p}^k / \mathfrak{p}^{\ell}}} \chi(xy) = \displaystyle\sum_{\substack{x \in \mathfrak{p}^{1 - \ell} / \mathfrak{p}^n - \mathfrak{p}^{1 - \ell + 1} / \mathfrak{p}^n \\ y \in \mathfrak{p}^k / \mathfrak{p}^{\ell}}} \chi(xy) + \displaystyle\sum_{\substack{x \in \mathfrak{p}^{1 - \ell+1} / \mathfrak{p}^n - \mathfrak{p}^{1 - \ell + 2} / \mathfrak{p}^n \\ y \in \mathfrak{p}^k / \mathfrak{p}^{\ell}}} \chi(xy)$$ $$ + \displaystyle\sum_{\substack{x \in \mathfrak{p}^{1 - \ell+2} / \mathfrak{p}^n - \mathfrak{p}^{1 - \ell + 3} / \mathfrak{p}^n \\ y \in \mathfrak{p}^k / \mathfrak{p}^{\ell}}} \chi(xy) + ... + \displaystyle\sum_{\substack{x \in \mathfrak{p}^n / \mathfrak{p}^n \\ y \in \mathfrak{p}^k / \mathfrak{p}^{\ell}}} \chi(xy)$$

Note that if $val(x) = i$, we require $y \in \mathfrak{p}^{-i} / \mathfrak{p}^{\ell}$ in order to force $xy \in \mathfrak{o}$.  Moreover, note that if $val(x) = i - \ell$, and $y$ ranges over $\mathfrak{p}^{\ell-i} / \mathfrak{p}^{\ell}$, then $xy$ ranges over $\mathfrak{o} / \mathfrak{p}^i$.  Coupling this with the fact that the sum of a nontrivial character over a group vanishes, we may compute the above sums easily.  We note that the evaluation of these sums can vary depending on $k, \ell, $ and $n$, as one can check.

We now assume $1 - \ell \leq m$.  Then there is no analogous term $A$ as above that we need to evaluate, and so the evaluation of the sum $$\displaystyle\sum_{\substack{x \in \mathfrak{p}^{m} / \mathfrak{p}^{n} \\ y \in \mathfrak{p}^k / \mathfrak{p}^{\ell}}} \chi(xy)$$ is analogous to the sum $B$ above.

Case ii) Assume $n + \ell < -2$.  In this case, it's never possible that $xy \in \mathfrak{o}$ unless $x$ or $y$ is zero.  Therefore, $$\displaystyle\sum_{\substack{x \in \mathfrak{p}^{m} / \mathfrak{p}^{n} \\ y \in \mathfrak{p}^k / \mathfrak{p}^{\ell}}} \chi(xy) = \displaystyle\sum_{\substack{x = 0 \\ y \in \mathfrak{p}^k / \mathfrak{p}^{\ell}}} \chi(xy) + \displaystyle\sum_{\substack{x \in \mathfrak{p}^{m} / \mathfrak{p}^{n} \\ y = 0}} \chi(xy) - \chi(0*0)$$  We subtracted $\chi(0*0)$ since we have double counted the term $\chi(0*0)$ in the right hand side of the equality. Thus, $$\displaystyle\sum_{\substack{x \in \mathfrak{p}^{m} / \mathfrak{p}^{n} \\ y \in \mathfrak{p}^k / \mathfrak{p}^{\ell}}} \chi(xy) = \left(\displaystyle\sum_{ y \in \mathfrak{p}^k / \mathfrak{p}^{\ell}} 1 \right) + \left(\displaystyle\sum_{x \in \mathfrak{p}^{m} / \mathfrak{p}^{n}} 1 \right) - 1 = q^{n-m} + q^{\ell-k} - 1$$  This finishes the case $n + \ell < -2$.

Case 5) Suppose $m + k < a, i = b$.  Then $$\displaystyle\sum_{ z \in \mathfrak{p}^i / \mathfrak{p}^j } \chi(\varpi^{-b} z) = 0$$ since $i = b$ (as in a previous case).  Moreover, in order to get $\varpi^{-a} xy + \varpi^{-b} z$ to be in $\mathfrak{o}$ we need the negative valuation terms of $\varpi^{-a} xy$ to be zero, since $\varpi^{-b} z$ is always in $\mathfrak{o}$.  Therefore, we are forced to take values of $x,y$ such that $\varpi^{-a} xy \in \mathfrak{o}$.  Over these values of $x,y$, we get $\chi(\varpi^{-a} xy + \varpi^{-b} z) = \chi(\varpi^{-a} xy) \chi(\varpi^{-b} z)$.  Therefore, $$\displaystyle\sum_{\substack{x \in \mathfrak{p}^m / \mathfrak{p}^n \\ y \in \mathfrak{p}^k / \mathfrak{p}^{\ell} \\ z \in \mathfrak{p}^i / \mathfrak{p}^j }} \chi(\varpi^{-a} xy + \varpi^{-b} z) = 0$$

Case 6) Suppose $m+ k < a, i < b$.  Recall from assumptions i) and ii) at the beginning of section \ref{gammasections}, we may assume that $n + \ell > a$ and $j > b$.  Since $m + k < a, i < b$, we will have negative $\varpi$ powers in both $\varpi^{-a} xy$ and $\varpi^{-b} z$.  Consider the negative valuation part of a term of the form $\varpi^{-a} xy + \varpi^{-b} z$.  We need this negative valuation part to be zero.  But once the negative valuation part of this is zero, we are free to let the rest of $x,y,z$ vary.  Fix possible negative valuation parts, denoted $(xy)_{-}$ and $z_{-}$, of $\varpi^{-a} xy$ and $\varpi^{-b} z$, respectively.  We compute the contribution to $$\displaystyle\sum_{\substack{x \in \mathfrak{p}^m / \mathfrak{p}^n \\ y \in \mathfrak{p}^k / \mathfrak{p}^{\ell} \\ z \in \mathfrak{p}^i / \mathfrak{p}^j }} \chi(\varpi^{-a} xy + \varpi^{-b} z)$$ of all terms $\chi(\varpi^{-a} xy + \varpi^{-b} z)$ such that $\varpi^{-a} xy$ has negative valuation $(xy)_{-}$ and $\varpi^{-b} z$ has negative valuation $(z)_{-}$, where $(xy)_{-} = -(z)_{-}$ (this last equality is forced upon us since otherwise $\chi$ will vanish).  Since $(xy)_{-} = -(z)_{-}$, this contribution is $$\displaystyle\sum_{\substack{x \in \mathfrak{p}^m / \mathfrak{p}^n, y \in \mathfrak{p}^k / \mathfrak{p}^{\ell}, z \in \mathfrak{p}^i / \mathfrak{p}^j \\ \mathrm{such \ that} \ xy \in \mathfrak{p}^a, z \in \mathfrak{p}^b }} \chi(\varpi^{-a} xy + \varpi^{-b} z) = \displaystyle\sum_{\substack{x \in \mathfrak{p}^m / \mathfrak{p}^n, y \in \mathfrak{p}^k / \mathfrak{p}^{\ell}, z \in \mathfrak{p}^i / \mathfrak{p}^j \\ \mathrm{such \ that} \ xy \in \mathfrak{p}^a, z \in \mathfrak{p}^b }} \chi(\varpi^{-a} xy) \chi(\varpi^{-b} z) = 0.$$
This argument holds regardless of $(xy)_{-}$ and $(z)_{-}$.  Thus, in the end, we are summing up a bunch of zeroes, so we finally get that $$\displaystyle\sum_{\substack{x \in \mathfrak{p}^m / \mathfrak{p}^n \\ y \in \mathfrak{p}^k / \mathfrak{p}^{\ell} \\ z \in \mathfrak{p}^i / \mathfrak{p}^j }} \chi(\varpi^{-a} xy + \varpi^{-b} z) = 0$$

\subsection{The calculation of $\Xi(x)$}

In this section we calculate the terms of the form $\Xi(x)$ in a much more general context.  In particular, we shall calculate $$\displaystyle\sum_{\substack{x \in \mathfrak{p}^m / \mathfrak{p}^n, y \in \mathfrak{p}^k / \mathfrak{p}^{\ell} \\ \mathrm{such \ that} \ xy \in \mathfrak{p}^c}} \chi(x) \chi(y)$$ We will calculate this sum in complete generality.  We split up the calculation of the above sum into various cases.

Case 1): Suppose $n + \ell \leq c + 1$.  Then it's never the case that $xy \in \mathfrak{p}^c$ unless at least one of $x,y$ are zero.   Therefore, if $n + \ell \leq c + 1$, we get
$$\displaystyle\sum_{\substack{x \in \mathfrak{p}^m / \mathfrak{p}^n, y \in \mathfrak{p}^k / \mathfrak{p}^{\ell} \\ \mathrm{such \ that} \ xy \in \mathfrak{p}^c}} \chi(x) \chi(y) = \displaystyle\sum_{y \in \mathfrak{p}^k / \mathfrak{p}^{\ell} \ } \chi(y) + \displaystyle\sum_{y \in \mathfrak{p}^m / \mathfrak{p}^n \ } \chi(x) - 1$$ which is easily calculatable, each sum being $0$ or a power of $q$, depending on $k,\ell,m,n$.  Note that we subtracted $\chi(0)\chi(0) = 1$ since this is $\chi(x) \chi(y)$ when $x = y = 0$ and we have double counted this term when we added the $x =0$ and $y = 0$ sums above.

Case 2): Now assume that $n + \ell > c+1$.  Assume furthermore that $1 - \ell + c > m$.  We will separate the sum into two parts.  We write
$$\displaystyle\sum_{\substack{x \in \mathfrak{p}^m / \mathfrak{p}^n, y \in \mathfrak{p}^k / \mathfrak{p}^{\ell} \\ \mathrm{such \ that} \ xy \in \mathfrak{p}^c}} \chi(x) \chi(y) =
\displaystyle\sum_{\substack{x \in \mathfrak{p}^{1-\ell+c} / \mathfrak{p}^n, y \in \mathfrak{p}^k / \mathfrak{p}^{\ell} \\ \mathrm{such \ that} \ xy \in \mathfrak{p}^c}} \chi(x) \chi(y) + $$
$$\displaystyle\sum_{\substack{x \in ( \mathfrak{p}^m / \mathfrak{p}^n - \mathfrak{p}^{1-\ell+c} / \mathfrak{p}^n ), y \in \mathfrak{p}^k / \mathfrak{p}^{\ell} \\ \mathrm{such \ that} \ xy \in \mathfrak{p}^c}} \chi(x) \chi(y).$$   (note that we can write $x \in ( \mathfrak{p}^m / \mathfrak{p}^n - \mathfrak{p}^{1-\ell+c} / \mathfrak{p}^n )$ since we assumed that $1 - \ell + c > m$).  We first compute the second sum on the right hand side.  In this sum, we will never have $xy \in \mathfrak{p}^c$ unless at least one of $x,y$ is zero.  But $x$ can't be zero, since $x \in \mathfrak{p}^m / \mathfrak{p}^n - \mathfrak{p}^{1-\ell+c} / \mathfrak{p}^n $.  Therefore, we must have $y = 0$.  Therefore, the second sum is $$\displaystyle\sum_{\substack{x \in ( \mathfrak{p}^m / \mathfrak{p}^n - \mathfrak{p}^{1-\ell+c} / \mathfrak{p}^n ), y \in \mathfrak{p}^k / \mathfrak{p}^{\ell} \\ \mathrm{such \ that} \ xy \in \mathfrak{p}^c}} \chi(x) \chi(y) = \displaystyle\sum_{x \in ( \mathfrak{p}^m / \mathfrak{p}^n - \mathfrak{p}^{1-\ell+c} / \mathfrak{p}^n )} \chi(x) \chi(0) = $$ $$\displaystyle\sum_{x \in ( \mathfrak{p}^m / \mathfrak{p}^n)} \chi(x)  - \displaystyle\sum_{x \in ( \mathfrak{p}^{1-\ell+c} / \mathfrak{p}^n) } \chi(x)$$ which is easy to calculate, each sum being $0$ or a power of $q$, depending on $m,n,\ell,c$.

So we now need to calculate the first sum $$\displaystyle\sum_{\substack{x \in \mathfrak{p}^{1-\ell+c} / \mathfrak{p}^n, y \in \mathfrak{p}^k / \mathfrak{p}^{\ell} \\ \mathrm{such \ that} \ xy \in \mathfrak{p}^c}} \chi(x) \chi(y)$$
We separate this into cases:

Case A: $1 - \ell + c > 0$.  We will calculate the sum by fixing the valuation of $x$.  Namely, $$\displaystyle\sum_{\substack{x \in \mathfrak{p}^{1-\ell+c} / \mathfrak{p}^n, y \in \mathfrak{p}^k / \mathfrak{p}^{\ell} \\ \mathrm{such \ that} \ xy \in \mathfrak{p}^c}} \chi(x) \chi(y) = \displaystyle\sum_{\substack{x \in (\mathfrak{p}^{1 - \ell + c} / \mathfrak{p}^n - \mathfrak{p}^{1 - \ell + c + 1} / \mathfrak{p}^n), y \in \mathfrak{p}^k / \mathfrak{p}^{\ell} \\ \mathrm{such \ that} \ xy \in \mathfrak{p}^c}} \chi(x) \chi(y) + $$ $$\displaystyle\sum_{\substack{x \in (\mathfrak{p}^{1 - \ell + c+1} / \mathfrak{p}^n - \mathfrak{p}^{1 - \ell + c + 2} / \mathfrak{p}^n) \\ y \in \mathfrak{p}^k / \mathfrak{p}^{\ell} \\ \mathrm{such \ that} \ xy \in \mathfrak{p}^c}} \chi(x) \chi(y) + $$ $$\displaystyle\sum_{\substack{x \in (\mathfrak{p}^{1 - \ell + c+2} / \mathfrak{p}^n - \mathfrak{p}^{1 - \ell + c + 3} / \mathfrak{p}^n) \\ y \in \mathfrak{p}^k / \mathfrak{p}^{\ell} \\ \mathrm{such \ that} \ xy \in \mathfrak{p}^c}} \chi(x) \chi(y) + ... + \displaystyle\sum_{\substack{x \in \mathfrak{p}^n / \mathfrak{p}^n, y \in \mathfrak{p}^k / \mathfrak{p}^{\ell} \\ \mathrm{such \ that} \ xy \in \mathfrak{p}^c}} \chi(x) \chi(y)$$

Note that if $val(x) = 1 - \ell + c + i$, then in order for $xy$ to be in $\mathfrak{p}^c$, we must have $y \in \mathfrak{p}^{\ell-i-1}$.  Coupling this with the fact that $\chi(x) = 1$ in every summation since we have assumed that $1 - \ell + c > 0$, we may compute the above sums easily.  We note that the evaluation of these sums can vary depending on $c, k, \ell, $ and $n$, as one can check.

\

Case B: Suppose $1 - \ell + c \leq 0$ and $\ell > 0$.  Suppose first that $n \leq 0$.  Then $\chi(x) \chi(y) \neq 0$ iff $x = 0$ since $\chi$ is zero on elements of negative valuation.  Therefore, if $n \leq 0$, $$\displaystyle\sum_{y \in \mathfrak{p}^k / \mathfrak{p}^{\ell}} \chi(0) \chi(y) = \displaystyle\sum_{y \in \mathfrak{p}^k / \mathfrak{p}^{\ell}} \chi(y)$$ which is easy to calculate.  If $n > 0$, then again, since $\chi$ is zero on elements of negative valuation, and since $1 - \ell + c \leq 0$, so we get $$\displaystyle\sum_{\substack{x \in \mathfrak{p}^{1-\ell+c} / \mathfrak{p}^n, y \in \mathfrak{p}^k / \mathfrak{p}^{\ell} \\ \mathrm{such \ that} \ xy \in \mathfrak{p}^c}} \chi(x) \chi(y) = \displaystyle\sum_{\substack{x \in \mathfrak{o} / \mathfrak{p}^n, y \in \mathfrak{p}^k / \mathfrak{p}^{\ell} \\ \mathrm{such \ that} \ xy \in \mathfrak{p}^c}} \chi(x) \chi(y)$$
In the case that $k \geq c$, $xy \in \mathfrak{p}^c$ always holds, so
$$\displaystyle\sum_{\substack{x \in \mathfrak{o} / \mathfrak{p}^n, y \in \mathfrak{p}^k / \mathfrak{p}^{\ell} \\ \mathrm{such \ that} \ xy \in \mathfrak{p}^c}} \chi(x) \chi(y) = \displaystyle\sum_{x \in \mathfrak{o} / \mathfrak{p}^n} \chi(x) \displaystyle\sum_{y \in \mathfrak{p}^k / \mathfrak{p}^{\ell}} \chi(y)$$ which is easy to calculate, depending on whether or not $n$ is zero.

Suppose that $k < c$.  We break the sum
$$\displaystyle\sum_{\substack{x \in \mathfrak{o} / \mathfrak{p}^n, y \in \mathfrak{p}^k / \mathfrak{p}^{\ell} \\ \mathrm{such \ that} \ xy \in \mathfrak{p}^c}} \chi(x) \chi(y)$$
into two sums :

$$\displaystyle\sum_{\substack{x \in \mathfrak{o} / \mathfrak{p}^n, y \in \mathfrak{p}^k / \mathfrak{p}^{\ell} \\ \mathrm{such \ that} \ xy \in \mathfrak{p}^c}} \chi(x) \chi(y) = \displaystyle\sum_{\substack{x \in \mathfrak{o} / \mathfrak{p}^n, y \in \mathfrak{p}^c / \mathfrak{p}^{\ell} \\ \mathrm{such \ that} \ xy \in \mathfrak{p}^c}} \chi(x) \chi(y) + \displaystyle\sum_{\substack{x \in \mathfrak{o} / \mathfrak{p}^n, y \in (\mathfrak{p}^k / \mathfrak{p}^{\ell} - \mathfrak{p}^c / \mathfrak{p}^{\ell}) \\ \mathrm{such \ that} \ xy \in \mathfrak{p}^c}} \chi(x) \chi(y)$$
Note that we assumed that $1 - \ell + c \leq 0$, so we can talk about $\mathfrak{p}^c / \mathfrak{p}^{\ell}$.  Moreover, since we assumed that $k < c$, we can talk about $y \in (\mathfrak{p}^k / \mathfrak{p}^{\ell} - \mathfrak{p}^c / \mathfrak{p}^{\ell})$.  Call the first sum D and the second sum E.  We first analyze E.  Notice that since are subtracting all elements $y$ of valuation $\geq c$, if we take an element $x$ of valuation zero, $xy$ can never be in $\mathfrak{p}^c$. Therefore, $$E = \displaystyle\sum_{\substack{x \in \mathfrak{p} / \mathfrak{p}^n, y \in (\mathfrak{p}^k / \mathfrak{p}^{\ell} - \mathfrak{p}^c / \mathfrak{p}^{\ell}) \\ \mathrm{such \ that} \ xy \in \mathfrak{p}^c}} \chi(x) \chi(y)$$ which equals
$$\displaystyle\sum_{\substack{x \in \mathfrak{p} / \mathfrak{p}^n, y \in \mathfrak{p}^k / \mathfrak{p}^{\ell} \\ \mathrm{such \ that} \ xy \in \mathfrak{p}^c}} \chi(x) \chi(y) - \displaystyle\sum_{\substack{x \in \mathfrak{p} / \mathfrak{p}^n, y \in \mathfrak{p}^c / \mathfrak{p}^{\ell} \\ \mathrm{such \ that} \ xy \in \mathfrak{p}^c}} \chi(x) \chi(y).$$ The first of these two sums can be handled by Case A, and the second sum equals $$\displaystyle\sum_{y \in \mathfrak{p}^c / \mathfrak{p}^{\ell}} \chi(y),$$ which is easy to calculate, and depends on $c$ and $\ell$.

We now handle the first sum D. It's always the case that $xy \in \mathfrak{p}^c$ in this sum, so we get $$\displaystyle\sum_{\substack{x \in \mathfrak{o} / \mathfrak{p}^n, y \in \mathfrak{p}^c / \mathfrak{p}^{\ell} \\ \mathrm{such \ that} \ xy \in \mathfrak{p}^c}} \chi(x) \chi(y) = \displaystyle\sum_{\substack{x \in \mathfrak{o} / \mathfrak{p}^n \\ y \in \mathfrak{p}^c / \mathfrak{p}^{\ell}}} \chi(x) \chi(y) = \displaystyle\sum_{x \in \mathfrak{o} / \mathfrak{p}^n} \chi(x) \displaystyle\sum_{y \in \mathfrak{p}^c / \mathfrak{p}^{\ell}} \chi(y) = 0$$  since $$\displaystyle\sum_{x \in \mathfrak{o} / \mathfrak{p}^n} \chi(x) = 0$$ (recall that we are in the case that $n > 0$). Thus, $D = 0$.

Case C: Suppose $1 - \ell + c \leq 0$ and $\ell \leq 0$.  Then since $\chi$ is zero on negative valuation terms, the $y$ terms don't contribute unless $y = 0$.  Then $$\displaystyle\sum_{\substack{x \in \mathfrak{p}^{1-\ell+c} / \mathfrak{p}^n, y \in \mathfrak{p}^k / \mathfrak{p}^{\ell} \\ \mathrm{such \ that} \ xy \in \mathfrak{p}^c}} \chi(x) \chi(y) = \displaystyle\sum_{x \in \mathfrak{p}^{1-\ell+c} / \mathfrak{p}^n} \chi(x) = \displaystyle\sum_{x \in \mathfrak{o} / \mathfrak{p}^n} \chi(x)$$ which is easy to calculate, depending on whether or not $n$ is zero.

\

Finally, the above analysis for Case 2) assumed that $1 - \ell + c > m$.  The case that $1 - \ell + c \leq m$ is simpler and similar to the case of $1 - \ell + c \leq m$, as one can check.  We note that there is no sum of the form $$\displaystyle\sum_{\substack{x \in ( \mathfrak{p}^m / \mathfrak{p}^n - \mathfrak{p}^{1-\ell+c} / \mathfrak{p}^n ), y \in \mathfrak{p}^k / \mathfrak{p}^{\ell} \\ \mathrm{such \ that} \ xy \in \mathfrak{p}^c}} \chi(x) \chi(y)$$ that we need to evaluate, in the case that $1 - \ell + c \leq m$.

\end{document}